# *The Reconstruction of Theaetetus' Theory of Ratios of Magnitudes*
## by Stelios Negrepontis and Dimitrios Protopapas

*Abstract*. In the present chapter, we obtain the reconstruction of Theaetetus' theory of ratios of magnitudes based, according to Aristotle's *Topics* 158b, on the definition of proportion in terms of equal anthyphairesis. Our reconstruction is built on the anthyphairetic interpretation of the notoriously difficult *Theaetetus* 147d6-e1 passage on Theaetetus' mathematical discovery of quadratic incommensurabilities, itself based on the traces it has left on Plato's philosophical definition of Knowledge in his dialogues *Theaetetus*, *Sophist* and *Meno*. Contrary to earlier reconstructions by Becker, van der Waerden and Knorr, our reconstruction reveals a theory that (a) applies only to the restricted class of pairs of magnitudes whose anthyphairesis is finite or eventually periodic, and (b) avoids the problematic use of Eudoxus' definition 4 of Book V of Euclid's *Elements*.

## **Section 1.** *Introduction*

Aristotle, in *Topics* 158b29-35, informs us of a theory of ratios of magnitudes based on an anthyphairetic definition of proportion that predates Eudoxus' theory. Since Becker (1933) this theory was attributed to Theaetetus (Section 2). The importance of this early theory cannot be overestimated, as it is the first step towards the construction of the real numbers, a project that was brought to fruition only in modern times by Dedekind in 1872, following closely Eudoxus' theory of ratios of magnitudes, the theory of ratios of magnitudes that replaced Theaetetus' theory. Attempts to a satisfactory understanding of the theory have been hampered by the realization that Eudoxus' condition (Book V, Definition 4 of the *Elements*) appears to be necessary for the reconstruction (Becker 1933, van der Waerden 1954, Knorr 1975). This apparent anachronism has led some researchers to believe either that Theaetetus' theory was not rigorous, treating the analogue of the crucial Proposition V.9 of the *Elements* as self-evident (Knorr), or even that such theory never really existed (Acerbi 2003, Saito 2003) (Section 3). Previous historians of Mathematics along with Platonic scholars have failed to realize the opportunities offered to modern readers of Plato for a double understanding of Theaetetus' Mathematics and Plato's Philosophy by the plan outlined in the *Theaetetus*. Plato sets the fundamental question on philosophical knowledge in the beginning of the dialogue, where (a) philosophical knowledge is correlated with Theaetetus' mathematical discovery of a general theorem of incommensurability, as reported in the *Theaetetus* 147d3-e1, urging in fact Theaetetus to approach the philosophical question on knowledge by imitating his mathematical discovery, and (b) philosophical knowledge is described (i) as True Opinion plus Logos or, equivalently, (ii) as Name plus Logos (Section 4). Some of these opportunities offered by Plato, but ignored by modern scholars, as philosophical imitations of Theaetetus' mathematical discovery consist in the divisions of the Angler, the Noble Sophistry and the Sophist in the *Sophist*, where each is described as providing knowledge by Name plus Logos, and in the *Meno*, which provides knowledge of the diameter to the side of a square, described as True Opinion plus Logos. Both the entities described as Plato's Name plus Logos in the *Sophist* and as True

Opinion plus Logos in the *Meno* reveal lucidly that Logos, the fundamental concept of Plato's philosophy, is an imitation of the Logos Criterion for anthyphairetic periodicity, namely the condition that the ratio of two successive anthyphairetic remainders is equal to the ratio of two Subsequent successive remainders (Section 2). We conclude that
(a) contrary to the arithmetical reconstruction of Theaetetus' proof regarding quadratic incommensurabilities provided by Euclid in Proposition X.9 of the *Elements*, which both van der Waerden (1954) and Knorr (1975) suggest (Section 9), Theaetetus' theorem reported in passage 147d3-e1 states that

> if lines a and b satisfy $a^2 = N \cdot b^2$ for non-square number N, or more generally $A \cdot a^2 = C \cdot b^2$ for natural numbers A, C, such that A·C is not square, then the anthyphairesis of a to b satisfies the Logos Criterion for eventual periodicity (and hence, by Proposition X.2, a and b are incommensurable) marking the transition from the Pythagorean method of finitization by the preservation of Gnomons to the Logos Criterion and the birth of Theaetetus' theory of ratios of magnitudes (Sections 6, 7);
(b) contrary to Burnyeat's (1978) claim (Section 9), Plato had a strong motive to exhibit the method of proof of Theaetetus, since Plato's Logos is the philosophical imitation of Theaetetus' Logos Criterion, revealing the anthyphairetic nature and structure of an intelligible Being, and thus forming the basis of Plato's own philosophy (Section 5);
(c) contrary to the insistent interpretation of Logos by Platonic scholars, as "account", "discourse", "rational definition","rational explanation" and "reason" (Section 9), we suggest that the meaning of Plato's Logos is the philosophical analogue of Theaetetus' Logos Criterion for anthyphairetic periodicity, resulting in geometric incommensurability and philosophical Knowledge and Oneness of the intelligible Beings (Section 5);
(d) contrary to claims made by Acerbi and Saito, Theaetetus' theory of proportion for magnitudes is certainly existent, since essential use of this theory is made in Theaetetus' theorem on quadratic incommensurabilities (Section 7), in the Theaetetean Book X of the *Elements*, including the palindromicity of anthyphairesis of a to b if $a^2 = Nb^2$ for a non-square number N and the solution to Pell's problem in relation to the Platonic dialogue *Statesman* (Negrepontis, Farmaki, Brokou, 2021) and through its imitation in Plato's philosophy in the *Sophist*, the *Meno* and other dialogues (Section 5); and,
(e) contrary to the reconstructions put forward by Becker (1933), van der Waerden (1954) and Knorr (1975) (Section 3), Theaetetus' theory of ratios of magnitudes is

> NOT about the class of all ratios of magnitudes, as Eudoxus' theory of ratios, presented in Book V of the *Elements*,
> BUT is strictly about the limited class of all ratios with finite or eventually periodic anthyphairesis (Section 8).

Once it is made clear that Theaetetus' theory of ratios is strictly and only for ratios with finite or eventually periodic anthyphairesis (certainly not because the definition of equal anthyphairesis does not allow it, but because ratios with non-periodic anthyphairesis are not wanted), then we will have no problem in reconstructing the theory, all the while avoiding the use of the anachronistic Eudoxus' condition of Definition 4 of Book V of the *Elements* (Section 10).

We expect to provide further support towards the confirmation of our reconstruction in a separate publication, by showing that the proofs of all the propositions in the Theaetetean Book X of the *Elements* do not need Eudoxus' theory of ratios which relies on Definition 4 of Book V of the *Elements*, but only the limited Theaetetus' theory of ratios applied solely to ratios with terms in eventually periodic anthyphairesis.

**Section 2**. *Aristotle, in* Topics *158b, reports of a pre-Eudoxean theory of ratios of magnitudes based upon equality of anthyphairesis, which, even though Aristotle does not mention its originator, is generally attributed to Theaetetus*

**2.1.** *A Pre-Eudoxean theory of Proportion for magnitudes based on the equality of anthyphairesis, is reported by Aristotle in the* Topics *158b29-35*

Book V of Euclid's *Elements* is devoted to the presentation of Eudoxus' theory of ratios of magnitudes, on which Books X, XII, and XIII of the Elements are based. However, Aristotle in *Topics* 158b29-35 reports of a clearly *preceding* theory of ratios of magnitudes, based on the definition of proportion in terms of equality of antanairesis/anthyphairesis.

> ἔοικε δὲ καὶ ἐν τοῖς μαθήμασιν
> ἔνια δι' ὁρισμοῦ ἔλλειψιν οὐ ῥᾳδίως γράφεσθαι,
> οἷον ὅτι ἡ παρὰ τὴν πλευρὰν τέμνουσα τὸ ἐπίπεδον
> ὁμοίως διαιρεῖ τήν τε γραμμὴν καὶ τὸ χωρίον.
> τοῦ δὲ ὁρισμοῦ ῥηθέντος εὐθέως φανερὸν τὸ λεγόμενον·
> τὴν γὰρ αὐτὴν ἀνταναίρεσιν ἔχει τὰ χωρία καὶ αἱ γραμμαί·
> ἔστι δ' ὁρισμὸς τοῦ αὐτοῦ λόγου οὗτος.
>
> It also appears that in mathematics,
> it is not easy to construct proofs of some things because of the deficiency of a definition,
> for instance that the line cutting a plane figure parallel to its side
> divides the area and the line similarly.
> But once the definition has been stated, the proposition is evident at once.
> For the areas and the lines have the same reciprocal subtraction
> (antanairesis, anthyphairesis):
> and this is the definition of 'same ratio'.

[[Aristotle *Topics* Books I and VIII, with excerpts from related texts. *Translated with a Commentary by* Robin Smith, Clarendon Press, Oxford, 1997]]

That the term 'antanairesis' in Aristotle's *Topics* coincides with the term 'anthyphairesis' in Euclid's *Elements* (Propositions VII.1 & 2, X.2) is confirmed by Alexander Aphrodisieus, *Commentary to Aristotle's Topics* 545, 16-17.

We can safely reconstruct the definition that Aristotle has in mind as follows.

Let a, b be magnitudes, a > b, with finite or infinite anthyphairesis

$a = k_0 \cdot b + c_1, c_1 < b,$
$b = k_1 \cdot c_1 + c_2, c_2 < c_1,$

$$c_1 = k_2 \cdot c_2 + c_3, \; c_3 < c_2,$$
$$\ldots$$
$$c_{n-1} = k_n \cdot c_n + c_{n+1}, \; c_{n+1} < c_n,$$
$$\ldots$$

We set

$$\text{Anth}(a, b) = [k_0, k_1, k_2, \ldots, k_n, \ldots].$$

Two dyads of magnitudes a, b, a > b and c, d, c > d have *equal anthyphairesis* if

$$\text{Anth}(a, b) = \text{Anth}(c, d).$$

The definition of proportion that Aristotle has in mind in his *Topics* passage is then the following:

*Definition.* If a, b, a > b, and c, d, c > d, are two dyads of magnitudes, we say that a/b = c/d if Anth(a, b) = Anth(c, d)

The Proposition that Aristotle has in mind is a simpler version of Proposition VI.1 of the *Elements*, and may be stated in modern notation as follows:

*Proposition.* If a, b, c are line segments, then a·c/b·c = a/b.

*Proof.* Let the anthyphairesis of a to b be the following:

$$a = k_0 \cdot b + c_1, \; c_1 < b,$$
$$b = k_1 \cdot c_1 + c_2, \; c_2 < c_1,$$
$$c_1 = k_2 \cdot c_2 + c_3, \; c_3 < c_2,$$
$$\ldots$$
$$c_{n-1} = k_n \cdot c_n + c_{n+1}, \; c_{n+1} < c_n,$$
$$\ldots$$

Then

$$a \cdot c = k_0 \cdot b \cdot c + c_1 \cdot c, \; c_1 \cdot c < b \cdot c,$$
$$b \cdot c = k_1 \cdot c_1 \cdot c + c_2 \cdot c, \; c_2 \cdot c < c_1 \cdot c,$$
$$c_1 \cdot c = k_2 \cdot c_2 \cdot c + c_3 \cdot c, \; c_3 \cdot c < c_2 \cdot c,$$
$$\ldots$$
$$c_{n-1} \cdot c = k_n \cdot c_n \cdot c + c_{n+1} \cdot c, \; c_{n+1} \cdot c < c_n \cdot c,$$
$$\ldots$$

It is clear that Anth(a, b) = Anth(a·c, b·c), hence a·c/b·c = a/b.

**2.2.** *Although Aristotle does not mention the originator, nevertheless the theory of proportion is normally attributed by modern historians of Greek Mathematics to Theaetetus, the great geometer of Plato's Academy*

Aristotle does not mention the discoverer of this theory of proportion of magnitudes. Becker, van der Waerden, and Knorr regard Theaetetus the originator of this theory of proportion:

Knorr, 1975,
>…an attractive thesis by O. Becker [1933]: that in the period before Eudoxus an attempt had been made to produce a general theory of proportions, valid for commensurable and incommensurable magnitudes alike, via a proportion-concept based on anthyphairesis. [Footnote 43 on p. 130] This thesis was also proposed by O. Toeplitz, H. Zeuthen and E. Dijksterhuis. For a discussion, see B. L. van der Waerden, *SA* 1954, pp. 175-179. p.125-126.
>…a proportion theory based on anthyphairesis was indeed the keystone of Theaetetus' theory of irrational lines. p.126.
>
>…we argue that Theaetetus introduced the anthyphairetic theory of proportion as the necessary auxiliary to his study of irrationals, not as a formal remedy of the 'crisis' in proportion theory occasioned by the discovery of incommensurability. p.290.

**Section 3**. *The problematic reliance of the reconstructions of Theaetetus' theory of ratios of magnitudes by Becker, 1933 and Knorr, 1975 on Eudoxus' condition (Definition V.4 of the* Elements*)*

Becker assumes that Theaetetus' theory of ratios of magnitudes, based on the definition of proportion in terms of equal anthyphairesis, is a theory for the same class of magnitudes as is Eudoxus' theory. Therefore, his reconstruction of Theaetetus' proof of the analogue of the crucial Proposition V.9 of the *Elements*, is inevitably in need of Eudoxus' condition (Definition V.4); van der Waerden essentially agrees with this reconstruction:

van der Waerden, 1954, p.177-178
>O. Becker has advanced an ingenious hypothesis for this. From the proportionality
>(1) a : b = c : d,
>one deduces first the equality of the areas
>(2) ad = bc
>then interchanges *b* and *c,* and finally returns to the proportionality
>(3) *a: c = b: d*
>It is this last step which involves the proposition on rectangles, that Aristotle talked about:
>(4) a : c = ad : cd,
>since
>a : c = ad : cd = cb : cd = b : d
>follows from (4) and (2).
>The fragments from Aristotle, which have been cited, make it very plausible,
>that (3) was indeed derived from (2) in this manner.
>The deduction of (2) from (1) involves another proposition. For, from (1)
>follows
>ad : bd = a : b = c : d = bc : bd.
>Now it remains to prove:
>P. If in a proportion the consequents are equal then the antecedents
>(ad and bc) are equal as well.
>The proof of P requires use of the so-called "lemma of Archimedes", that is

> formulated as follows: *Q. If A and B are comparable magnitudes (e.g. both line segments), and if A is less than B, then a certain multiple nA exceeds B.*

Note that what van der Waerden calls "lemma of Archimedes" is Definition 4 of Book V of the *Elements,* and we refer to it as "Eudoxus' condition".

van der Waerden, 1954, p.178-179
> All becomes clear now. Evidently, Theaetetus began his book with an exposition of the theory of proportions, based on the antanairesis-definition.
> Following his usual procedure, he started with lemmas which would be needed later on; among these is R = X 1. In propositions X 2, 3 he established the theory of the infinite or finite antanairesis, thus obtaining at the same time a criterium for the commensurability of two line segments or two areas.
> It is probable that the next thing was the theory of proportions, based on the antanairesis-definition to which Aristotle refers. Euclid omitted this part, because he had already given another theory of proportions in Book V.

Knorr, 1975, although not fully accepting Becker's reconstruction, nevertheless does not question Becker's assumption that Theaetetus' theory is about the same class of magnitudes as is the theory of Eudoxus, and that therefore a rigorous proof of the analogue of Proposition V.9 would need Eudoxus' condition.

> The simple elegance of the Euclidean proof in contrast to the cumbersome anthyphairetic version scarcely needs comment. Becker and the adherents to his thesis on the pre-Eudoxean anthyphairetic theory wish to maintain that this theorem was proved within that theory and made necessary for the first time the use of the convergence theory. I do not fully accept this view. For the algebra [of the anthyphairetic proof] is messy, inevitably so; but by Greek methods of algebra it would be horrendous, and I have no doubt in this case prohibitively difficult to effect formally. I find it significant that Theorem 6 [if for magnitudes A, B, C, we have A:C=B:C, then A=B (Elements V.9)] is absent from Book VII. This leads me to believe that in the initial stages in the theory of proportions, both for integers and for magnitudes (via anthyphairesis), this deceptively easy-looking step was omitted.
> But in the course of the complete formulation of the theory, and probably in context of Theorem 7 [Given four lines A, B, C, D; A:B=C:D if and only if AD=BC (*Elements* VII.19 and VI.16)] (as Becker maintains), the necessity of Theorem 6 was detected. The anthyphairetic verification must certainly have been abandoned, however, after it was discovered that certain principles worked out after an attempt at a proof-i.e. X.1 and V, Def. 4-could be turned around into an alternative definition of proportion (V, Def.5). Thus, the stage was reached that at which the whole theory of proportions could now be put into Eudoxean-Euclidean form.

Therefore, if we assume that Theaetetus' theory was about essentially all ratios of magnitudes, we are forced to accept either an anachronistic use of Eudoxus' condition or some sloppy, non-rigorous approach by Theaetetus. Both are highly implausible.

Thus, the previous reconstructions of Theaetetus' theory have a serious problem and, consequently, a reconstruction that is rigorous and avoids any use of Eudoxus' condition will have a serious advantage.
*A satisfactory reconstruction of Theaetetus' theory of ratios must not make any use of Eudoxus' condition (V, Definition 4).*

**Section 4.** *Our strategy to examine the central role that the imitation of Theaetetus' theorem, reported in the* Theaetetus *147d-148b, plays in Plato's philosophy, proves fruitful towards the understanding of both the method of Theaetetus' proof and the nature of Plato's philosophy*

**4.1.** *Theaetetus' account of the close connection between Theaetetus' mathematical discovery and Plato's philosophy*

We have no ancient mathematical sources at our disposal, informing us on the method used by Theaetetus to prove his quadratic incommensurabilities. On the other hand, we have Plato's dialogues intact and, as it has been shown by one of us (Negrepontis 2000, 2012, 2018), Plato's dialogues contain traces of these Mathematics at the very core of their philosophy. In particular, in his trilogy *Theaetetus-Sophist-Statesman*, Plato guides us to the connection between his philosophical theory of knowledge of intelligible Beings and Theaetetus' mathematical method of proving quadratic incommensurabilities.

*Step 1. The philosophical problem on knowledge,* Theaetetus *145d-147c*
In the beginning of the dialogue Theaetetus, Socrates sets the basic unifying problem that runs through the whole trilogy: How can we obtain the knowledge of the fundamental entities in Plato's philosophy, the intelligible Beings, also referred to as Platonic Ideas?

*Step 2. Theaetetus' realization that his recent mathematical discovery is similar to the philosophical problem on knowledge,* Theaetetus *147c-e*
After some inconclusive discussion, Theaetetus recognizes that Socrates' philosophical question appears to be related to a mathematical discovery that he recently made (with a friend of his). In *Theaetetus* 147d-e, Theaetetus describes his discovery but Plato's language is obscure and in need of careful and reasoned interpretation. Among earlier historians of Mathematics, including van der Waerden (1954) and Knorr (1975) and Platonic scholars there is general agreement:
(a) that Theodorus proved that if $a^2 = Nb^2$ for a non-square number N, $3 \leq N \leq 17$ (or, possibly, $3 \leq N < 17$), then a, b are incommensurable lines, and
(b) that Theaetetus discovered and proved a general incommensurability theorem, whose minimal content is the following:

> **Theaetetus'** Theorem (weak form).
> If a, b are line segments, satisfying $a^2 = Nb^2$ with N a non-square number, or more generally satisfying $a^2/b^2 = C/A$ with A, C natural numbers such that AC is a non-square number, then a, b are lines incommensurable/irrational to each other.

The exact statement and the method of proof of the theorem remains to be determined.

*Step 3. It is urged by Socrates that philosophical knowledge can be achieved by imitating Theaetetus' mathematical discovery,* Theaetetus *148c9-d7*
Shortly after, at *Theaetetus* 148c-d, Socrates exhorts Theaetetus to "try to imitate" his mathematical discovery in order to answer the philosophical question about the knowledge of the Ideas.

> ΣΩ. Θάρρει τοίνυν περὶ σαυτῷ καὶ τὶ οἴου Θεόδωρον λέγειν,
> προθυμήθητι δὲ παντὶ τρόπῳ
> τῶν τε ἄλλων πέρι καὶ ἐπιστήμης
> λαβεῖν λόγον τί ποτε τυγχάνει ὄν.?
> ΘΕΑΙ. Προθυμίας μὲν ἕνεκα, ὦ Σώκρατες, φανεῖται. 148c9-d3
>
> ΣΩ. Ἴθι δή—καλῶς γὰρ ἄρτι ὑφηγήσω—
> ***πειρῶ μιμούμενος*** (148d4-5) τὴν περὶ τῶν δυνάμεων ἀπόκρισιν,
> *ὥσπερ ταύτας πολλὰς οὔσας ἑνὶ εἴδει περιέλαβες,*
> *οὕτω καὶ τὰς πολλὰς ἐπιστήμας ἑνὶ λόγῳ προσειπεῖν*. 148d4-7

> Socrates. You may be reassured then, about Theodorus' account of you,
> and set your mind on finding a definition of knowledge, as of anything else,
> with all the zeal at your command.
> Theaetetus. If it depends on my zeal, Socrates, the truth will come to light.
>
> Socrates. Forward, then. On the way you have just shown so well.
> ***Take as a model*** your answer about the roots:
> just as you found a single character to embrace all that multitude,
> so now try to find a single formula that applies to the many kinds of knowledge.
> [translation by Cornford, 1957]

> socrates: So be confident about yourself, believe what Theodorus said about you,
> and commit yourself completely to getting an account of knowledge: what, exactly, is it?
> theaetetus: If commitment is what counts, Socrates, the answer will appear.
>
> socrates: Come on then. You've just given us a good start.
> ***Try mimicking*** your answer about powers:
> just as there were many of them, and yet you covered them all with a single form,
> try now to apply a single account to the many knowledges there are.
> [translation by Christopher Rowe, 2015]

*Step 4. It is further stated by Socrates that philosophical knowledge consists of Name plus Logos, equivalently of True Opinion plus Logos,* Theaetetus *201d-202c*
Later in the *Theaetetus* 201d-202c two equivalent definitions are given for the Knowledge of an intelligible Being (the first described with greater clarity in the *Sophist* 218c4-5):

> δεῖ δὲ ἀεὶ παντὸς πέρι τὸ πρᾶγμα αὐτὸ μᾶλλον διὰ λόγων ἢ τοὔνομα μόνον συνωμολογῆσθαι χωρὶς λόγου. Sophist 218c4-5

however, we ought always in every instance to come to agreement about the thing itself by "Logos"
rather than about the mere name without "Logos".
[translation by Fowler (1921), modifications by the authors]

*Definition* 1. The Knowledge (episteme) of an intelligible Being is not Name only but Logos, as well.

Equivalently,

ὅταν μὲν οὖν ἄνευ λόγου τὴν ἀληθῆ δόξαν τινός τις λάβῃ,
ἀληθεύειν μὲν αὐτοῦ τὴν ψυχὴν περὶ αὐτό, γιγνώσκειν δ' οὔ·
τὸν γὰρ μὴ δυνάμενον δοῦναί τε καὶ δέξασθαι λόγον ἀνεπιστήμονα εἶναι περὶ τούτου·
προσλαβόντα δὲ λόγον
δυνατόν τε ταῦτα πάντα γεγονέναι καὶ τελείως πρὸς ἐπιστήμην
ἔχειν *Theaetetus* 202b8-c5
When therefore a man acquires without "Logos" the true opinion about anything,
his mind has the truth about it, but has no knowledge;
for he who cannot give and receive the "Logos" of a thing is without knowledge of it;
but when he has acquired also the "Logos"
he may possibly have become all that I have said and may now be perfect in knowledge.
[translation by H.N. Fowler, modifications by the authors]

*Definition 2*. The Knowledge of an intelligible Being is True Opinion about this true Being plus Logos.

The attempts to understand these definitions of intelligible knowledge in the *Theaetetus* are not successful. Schematically we have the following

Interplay between Theaetetus' Mathematics and Plato's Philosophy in Plato's *Theaetetus*

| *philosophical question on Knowledge by Socrates* | *geometrical/mathematical discovery by Theaetetus* |
|---|---|
| *Step 1* 145d4-e9 Question set | |
| 145e-147c Inconsequential Discussion on the matter | |
| 147c4 Hint that knowledge is infinite | |
| | *Step 2* 147c7-d1 Realization by Theaetetus of a connection between the question on philosophical Knowledge and his mathematical discovery |
| | 147d3-148b2 |

|  | Account of Theaetetus' discovery |
|---|---|
| *Step 3* 148c9-d7<br>"peiro mimoumenos"<br>The urge to imitate Theaetetus' discovery<br>in order to answer the philosiophical question |  |
| *Step 4* 201d-202c<br>Knowledge<br>as Name plus Logos,<br>equivalently as True Opinion plus Logos |  |
| 206c-210d [end of dialogue]<br>Despite attempts, failure to understand<br>the true meaning of Knowledge and Logos |  |

**4.2.** *Outline of our strategy on the basis of Theaetetus' account*

Combining Steps 3 and 4, we conclude that philosophical knowledge
on the one hand is described as Name plus Logos, equivalently as True Opinion plus Logos,
on the other hand is described as imitation of Theaetetus' mathematical discovery.
The philosophical knowledge of the line a whose square doubles the area of the square with side b is given in the *Meno* as True Opinion plus Logos. In Section 5.1, we find that the philosophical terms True Opinion, Logos and Knowledge have a quite specific geometrical meaning:
*True Opinion* consists in the first two steps of the anthyphairesis of the diameter a to the side b of a square,
*Logos* consists in the Logos Criterion for the periodicity of the anthyphairesis of the diameter a to the side b of a square, and
*Knowledge* of the diameter a to the side b consists in the infinite periodic anthyphairesis with sequence of successive quotients the sequence [1, 2, 2, 2, …].
The philosophical knowledge of each of the intelligible Beings Angler, Noble Sophistry, and Sophist, is given in the *Sophist*, the sequel of the *Theaetetus* in the Platonic trilogy *Theaetetus-Sophist-Statesman* devoted to knowledge, as Name plus Logos. In Section 5.2, we find that the philosophical terms Name, Logos and Knowledge have a quite specific geometrical meaning:
*Knowledge* of the Angler consists in its definition in terms of a philosophical imitation of an anthyphairesis,
whose anthyphairetic remainders are the *Names,* including the Angler itself, and
which possesses *Logos*, the philosophical imitation of the Theaetetus Logos Criterion for periodic anthyphairesis.
Similar definitions are given for the Noble Sophistry and especially for the Sophist.
Since philosophical knowledge has also been described as the imitation of the theorem of Theaetetus, we expect that in fact Theaetetus proved not only the incommensurability of each, but in fact its anthyphairetic periodicity.

In Section 5.3 we arrive at the same conclusion by carefully examining the philosophical wording of the description of Theaetetus' mathematical discovery in the *Theaetetus* 147d3-e1 passage, and arguing that "the powers" in the 147d6-e1 passage should be rendered as a distributive plural, an interpretation that suggests that the method of proof of the incommensurabilities of the powers by Theodorus and Theaetetus was anthyphairetic. In particular, the philosophical statement that Theaetetus proved that the infinite multitude of parts in which each power is divided by anthyphairesis can be collected into One strongly suggests that every power is a self-similar One, a Oneness that is conceivable only under periodic anthyphairesis. We are thus led again to the conclusion that Theaetetus proved that the anthyphairesis of every power with respect to the one foot line is periodic, or more generally

> *Theaetetus' Theorem (strong form).*
> If a, b are line segments, satisfying $a^2 = Nb^2$ with N non-square number, or more generally satisfying $a^2/b^2 = C/A$ with A, C natural numbers such that AC is a non-square number,
> then the anthyphairesis of a to b is eventually periodic,
> hence a, b are lines incommensurable/irrational to each other.

In conclusion our approach proves to be quite fruitful. It yields both
(a) that Theaetetus' theorem reported in the *Theaetetus* 147d3-b2 states that every quadratic equation with arithmetical coefficients possesses an eventually periodic anthyphairesis, established in terms of the Theaetetus' Logos Criterion, and
(b) that Plato's intelligible Being/Idea is the philosophical analogue of a dyad of kinds, satisfying *Logos*, the philosophical analogue of Theaetetus' Logos Criterion for anthyphairetic periodicity.
In addition, it makes clear
(c) that the role of Theaetetus' theory of ratios for magnitudes is fundamental, both for Theaetetus' theorem, proved with an appeal to the Logos Criterion for anthyphairetic periodicity, and for Plato's account of philosophical Knowledge in terms of Logos, the philosophical imitation of Theaetetus' Logos Criterion.

**Section 5.** *Plato's imitation of Theaetetus' discovery in the* Sophist *reveals*
*(a) the anthyphairetic nature of Plato's philosophy, in particular the crucial concept of Logos as the philosophical imitation of the Logos Criterion for anthyphairetic periodicity, and*
*(b) the anthyphairetic method of Theaetetus' proof of quadratic incommensurabilities, in particular its role for the birth of Theaetetus' theory of proportion for magnitudes*

**5.1.** *The philosophic imitation of Theaetetus' theorem in Plato's* Sophist *218b5-221c4*

The anthyphairetic nature, of the Divisions in the *Sophist* and in the *Statesman*, in particular of the Logos in these Divisions, have been examined elsewhere (Negrepontis, 2000, 2012, 2018, 2023; Negrepontis, Farmaki, Brokou, 2021), and the reader is referred there for the detailed arguments. In the present section we go through the simpler Division in the *Sophist*, that of the Angler.

**5.1.1.** *The Knowledge of the intelligible Being Angler by Name and Logos of the intelligible Being Angler in the* Sophist *218b-221c*

The method of 'Name and Logos' (cf. *Theaetetus* 201e2-202b5, *Sophist* 218c1-5, 221a7-b2,

268c5-d5), is exemplified at the start of the *Sophist* 218b-221c by the definition of the Being Angler. The Angler is certainly a lowly paradigm of an intelligible Being, but even for lowly entities, there is an intelligible Being, cf. *Parmenides* 130b7-e3. In the scheme below, we reproduce the binary division process which leads to the Angler.

The Division of the Angler (*Sophist* 218b-221c)

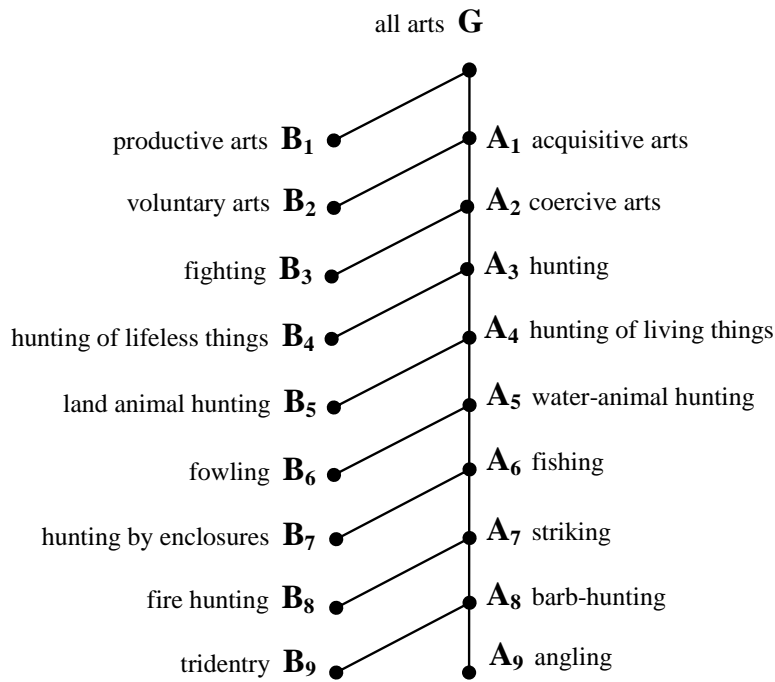

### 5.1.2. *The Logos Criterion of the Angler*

The description of the Logos Criterion of the Angler is contained in the *Sophist* 220e2-221c3 passage (which we have divided into two parts [A] and [B] for the sake of convenience):

> [A] '*Stranger* Then of striking which belongs to barb-hunting [A8],
> that part which proceeds *downward from above* ('anothen eis to kato'),
> is called, because tridents are chiefly used in it, tridentry [B9],
> I suppose….
> *Stranger* The kind that is characterized by the opposite sort of blow,
> which is practiced with a hook and strikes, ...
> and proceeds *from below upwards* ('katothen eis tounantion ano'),
> being pulled up by twigs and rods.
> By what name, Theaetetus, shall we say this ought to be called?
> *Theaetetus* I think our search is now ended and we have found the very thing we set before us a while ago as necessary to find.
> *Stranger*: Now, then, you and I are only agreed about the name of angling, [A9]'

220e2-221b1

> [B] 'but *we have acquired also a satisfactory 'Logos' of the thing itself.
> *For* (gar)
> of art as a whole, half was acquisitive,
> and of the acquisitive, half was coercive,
> and of the coercive, half was hunting,
> and of hunting, half was animal hunting,
> and of animal hunting, half was water hunting,
> and of water hunting [A5]
> *the whole part from below (to katothen tmema holon)* was fishing, [A6]
> and of fishing, half was striking,
> and of striking half was barb-hunting, [A8]
> and of this [A8]
> the part in which the blow is pulled
> from below upwards *(to peri ten katothen ano)* was angling. [A9] (221b1-c3).

(translation by H. N. Fowler 1921), with modifications by the author).

In [A], the opposing relation of Tridentry to Angling is carefully explained: all Fishing with a hook is divided into

> Tridentry (=Fishing with a trident), which is described as
> Fishing with a hook with an art that proceeds *from above downwards,* and
> Angling (=Fishing with a rod), which is described as
> Fishing with a hook with an art that proceeds *from below upwards.*
> We have pushed Division all the way to the Angling; thus, we have certainly found 'the name' of Angling.

But now, in [B], it is claimed that 'the Logos' of the Angling has been found, too.
The justification — the proof — that we have indeed found the Logos, too, is contained in the remainder [B], since this remaining part [B] starts with a 'for' ('gar'), and this justification can be seen to consist of:
(i) an accurate recounting of all the division steps, abbreviated in the sense that of the two species into which each genus is divided, only the one that contains the Angler is mentioned, while the species opposite to it is omitted;
(ii) a reminder that the last species, the angling, is characterized as the part of its genus that proceeds 'from below upwards'; and,
(iii) the ONLY new information (since (i) and (ii) are repetitions of things already contained in the Division and in [A]), which concerns the species of fishing, three steps before angling, and which informs us for the first time that this species is 'the whole part from below' of its genus.
Since this is an abbreviated account, there is no explicit information on the species opposite to 'fishing' — namely 'fowling' — but since 'fishing' was described not simply as 'the part from below' of its genus, but emphatically as 'the whole part from below', it follows that the opposite species — 'fowling' — must be characterized as 'the (whole) part from above' of the same genus. In fact, there can be no other justification for the presence of the term 'whole' in the description of 'fishing' with a view to arguing that we have obtained 'Logos', except to indicate and imply this description for its opposite species, 'fowling'.

We recall that the part [B] from the word 'for' ('gar') is explicitly a justification of the claim that we have succeeded in finding the 'Logos' of the Angling. We may then ask: what is the 'Logos' of the Angler that reasonably results from such a justification? There can really be only one answer: the 'Logos' we are looking for is the equality of the 'philosophic ratio' of Tridentry to Angling, namely the equality of the ratio 'of from above downwards to from below upwards', to the ratio 'of Fowling to Fishing'.

Thus:

> fowling B6/fishing A6 = tridentry B9/angling A9 =
> from above downwards/from below upwards.

Since the species Tridentry and Angling form a pair of opposite species, and the species Fowling and Fishing form another pair of opposite species in the Division Scheme for the Angler, the resulting 'Logos' bears a most uncanny similarity to the Logos Criterion for the periodicity of the anthyphairesis of geometric magnitudes, and of geometric powers in particular.

**5.1.3**. *The Name and Logos of the Angler*

The Name and Logos of the Angler thus takes the following form:

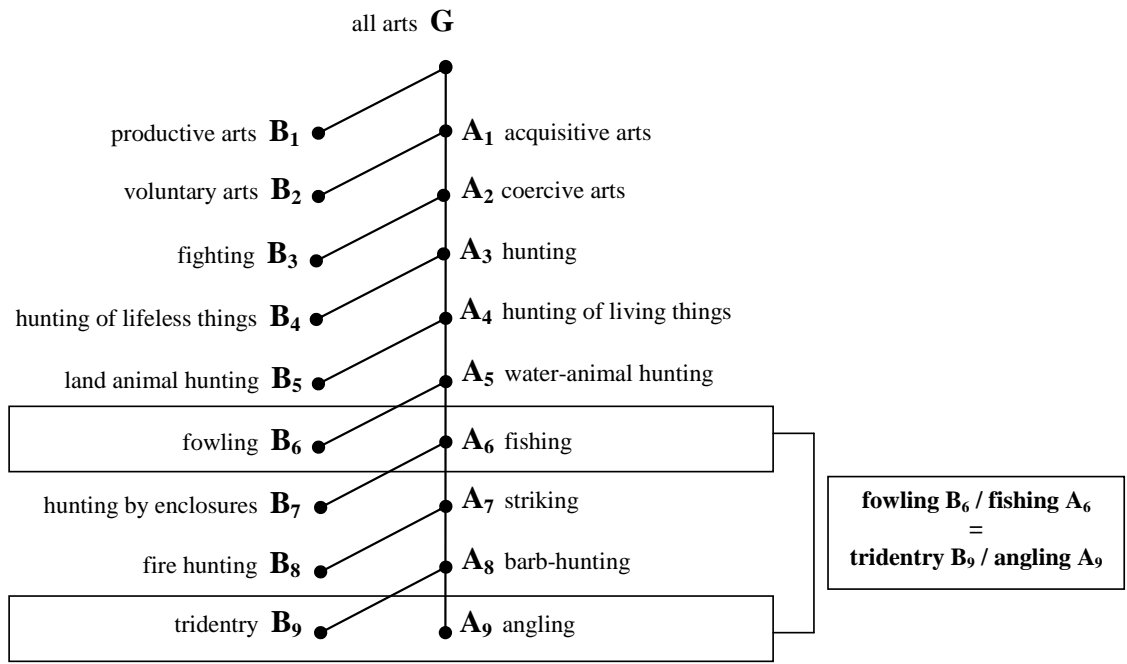

Thus the Name and Logos of the Angler consists of the anthyphairetic division together with the Logos, analogous to the Logos Criterion for periodicity for geometric anthyphairesis.

**5.2**. *Knowledge in the* Meno *80d-86a and 97a-98b as Recollection is the knowledge of the anthyphairesis of the diameter to the side of a square and is* *achieved by Theaetetus' Logos Criterion*

The nature of Plato's intelligible Being as a philosophic imitation of a dyad in periodic anthyphairesis is revealed in other Platonic dialogues, as well, such as the *Meno*, the *Republic*, the *Parmenides*, the *Timaeus* (cf. Negrepontis 2023, 2018, TA, Negrepontis & Kalisperi 2023, respectively). In the present section we outline the Knowledge of the diameter to the side of a square as obtained, in the *Meno* 80d-86e and 97a-98b, by True Opinion (the first two steps of the anthyphairesis of the diameter to the side) plus the Logos Criterion (instituting the periodicity of the anthyphairesis from there on).

**5.2.1.** *Proposition.* If $a \cdot d = b \cdot c$, then $Anth(a, b) = Anth(c, d)$.

*Proof.* Let $Anth(a, b) = [k_0, k_1, \ldots]$. We proceed by induction.
Thus $a = k_0 \cdot b + e_1$, with $e_1 < b$.
Then $c \cdot b = a \cdot d = k_0 \cdot b \cdot d + e_1 \cdot d$.
By Proposition I.44 of the *Elements*, there is a line segment $f_1$, such that $e_1 \cdot d = b \cdot f_1$.
Since $e_1 < b$, it is clear that $f_1 < d$.
Thus $c \cdot b = k_0 \cdot b \cdot d + e_1 \cdot d = k_0 \cdot b \cdot d + b \cdot f_1$; then $c = I_0 \cdot d + f_1$ with $f_1 < d$.
Thus $Anth(a, b) = [k_0, k_1, \ldots] = [k_0, Anth(b, e_1)]$, while $Anth(c, d) = [I_0, Anth(d, f_1)]$.
Thus, we have that $e_1 \cdot d = b \cdot f_1$ and we must prove that $Anth(b, e_1) = Anth(d, f_1)$.
We continue as in the first step, and we finish the proof by induction.

**5.2.2**. *Proposition (Pythagorean)*.
If b is a line and a is a line such that $a^2 = 2 \cdot b^2$,
then a is the diameter of the square with side b, and
$Anth(a, b) = [1, period(2)]$.

*Proof (Theaetetus).*
[Part I]
By the (isosceles) Pythagorean theorem, line a is constructed as the diameter of the square with given side b. Clearly $a > b$.
Since $a^2 = 2 \cdot b^2$, it follows that $a < 2 \cdot b$
(if we had $a \geq 2 \cdot b$, then $a^2 \geq 4 \cdot b^2 > 2 \cdot b^2 = a^2$, contradiction).
Thus $b < a < 2 \cdot b$. Hence,

$\quad a = b + c_1, c_1 < b.$

Claim. $b > 2 \cdot c_1, 3 \cdot b > 2 \cdot a$
(if we had $b \leq 2 \cdot c_1 = 2 \cdot a - 2 \cdot b$, then $3 \cdot b \leq 2 \cdot a$, hence $9 \cdot b^2 \leq 4 \cdot a^2 = 8 \cdot b^2$, contradiction). Hence,

$\quad b = 2 \cdot c_1 + c_2.$

$b > 2 \cdot (a - b), 3 \cdot b > 2 \cdot a.$

[Part II]
Then $c_1 = a - b$, hence $c_2 = b - 2 \cdot c_1 = b - 2 \cdot (a - b) = 3 \cdot b - 2 \cdot a$, hence

$b \cdot c_2 = b \cdot (3 \cdot b - 2 \cdot a) = 3 \cdot b^2 - 2 \cdot a \cdot b$, and
$c_1^2 = (a - b) \cdot (a - b) = a^2 + b^2 - 2 \cdot a \cdot b = 3 \cdot b^2 - 2 \cdot a \cdot b$.
Hence $b \cdot c_2 = c_1^2$. Since $b > c_1$, and

$$b \cdot c_2 = c_1^2,$$

clearly $c_2 < c_1$,
thus the second step is indeed anthyphairetic

$$b = 2 \cdot c_1 + c_2, \; c_2 < c_1.$$

By $b \cdot c_2 = c_1^2$ and Proposition 4.2.1, we have
Anth $(b, c_1)$ = Anth $(c_1, c_2)$, namely $b/c_1 = c_1/c_2$, according to the Theaetetean definition.
Then Anth $(c_1, c_2)$ = Anth $(b, c_1)$ = [2, Anth $(c_1, c_2)$]. Hence,
Anth $(a, b)$ =
[1, 2, Anth $(c_1, c_2)$] =
[1, 2, Anth $(b, c_1)$] =
[1, 2, 2, Anth $(c_1, c_2)$] =
… =
[1, period (2)].

**5.2.3**. *According to the* Meno, *the knowledge of the diameter to the side of a square is achieved by the knowledge of their infinite anthyphairesis*

By Proposition 5.2.2 we obtain *the complete knowledge* of the diameter a with respect to the side b of the square. Since the anthyphairesis of a to b is infinite, by Proposition X.2, a and b are incommensurable.
The computation depends on the representation of the first quotient $c_1$ as a – b and the second quotient $c_2$ as equal to 3b – 2a. These representations use the *side and diameter numbers* $q_n$, $p_n$; thus $c_1 = p_1 \cdot a - q_1 \cdot b = a - b$, $c_2 = q_2 \cdot b - p_2 \cdot a = 3 \cdot b - 2 \cdot a$.

**Section 6.** *The anthyphairetic reconstruction of Theaetetus' proof of quadratic incommensurabilities (*Theaetetus *147d3-e1), Part I: the computation of every step of the anthyphairesis with use of the Pythagorean Application of Areas, constructed by anthyphairetic substitution*

Plato's method of the acquisition of Knowledge of intelligible Beings,
described as True Opinion & Logos in the *Meno*,
described equivalently as Name & Logos in the *Sophist* and the *Statesman,* or as Division and Collection,
    according to our interpretation (outlined in Section 5),
    is the philosophical imitation of
    an anthyphairetic division followed by Logos, having the form e.g.
    for the Angler B6/A6 = B9/A9, for the diameter a to the side b, $b/c_1 = c_1/c_2$.

But Plato's method is also described in the *Theaetetus* 148c-d as the philosophical imitation of
Theaetetus' mathematical method of proof of his theorem on quadratic incommensurabilities.

Putting the two together, we certainly expect that Theaetetus' method of proof of his theorem on quadratic incommensurabilities will be anthyphairetic, in fact by anthyphairetic substitution. We present our anthyphairetic reconstruction of the proof of Theaetetus' theorem in the present Section 6 and in the next Section 7 (cf. also Negrepontis, Farmaki, Kalisperi, TA).

In Section 6 starting with an arithmetical quadratic equation, we present the construction of the anthyphairetic process, employing essentially the Pythagorean tools of the Application of Areas, the same tools that, as we have argued in the paper just mentioned, were used for the proof of the incommensurabilities by the Pythagoreans and by Theodorus. It is a process technically demanding, but not too different from the special cases worked out by Theodorus. The great novel idea that completes the proof will be given in the next Section 7.

### 6.1. *The Pythagorean Application of Areas in Book II of the* Elements

Emblematic within Book II of the *Elements* is Proposition II.4.
It has long been understood that the proper/original place of the Pythagorean theorem is not as Proposition I.47, proved with the help of Propositions I.37 and 1.41, but as Proposition II.4/5, proved with the help of Proposition II.4. This is strongly suggested by a variety of reasons, but mainly by the proof of the isosceles Pythagorean theorem described in Plato's *Meno* 82a-85b.
Furthermore, the method of Application of Areas appears in the *Elements* in general form, as Propositions VI.27 & 28 (in defect) and VI.29 (in excess), where use of Eudoxus' theory of ratios of magnitudes is made. On the other hand, Eudoxus' theory is presented in Book V and it has also long been understood, that the proper place of the original Pythagorean method, with no use of ratios, belongs to Book II of the *Elements*.
The two Pythagorean Propositions on Applications of Areas, in defect & in excess, are proved using the Pythagorean theorem and the Propositions II.5 & 6, respectively, which are variations of Proposition II.4, precisely suitable for these proofs. The Pythagorean Application of Areas in defect takes the position II.5/6 and in excess the position II.6/7.
Thus, the first part of Book II of the *Elements*, partially completed with the addition of the Pythagorean theorem and the Pythagorean Application of Areas has the following content:

*Propositions II.1, 2, 3.* If $a_1, a_2, ..., a_n$ and b are lines, then
$(a_1 + a_2 + ... + a_n) \cdot b = a_1 \cdot b + a_2 \cdot b + ... + a_n \cdot b$.
The distributive rule of the "multiplication" of magnitudes.

*Proposition II.4.*
If a, b are lines, then $(a + b)^2 = a^2 + b \cdot (b + 2a)$.
[Proof with the use of Gnomons.]
Proposition II.4 is the fundamental Proposition of the Pythagorean Book II of the *Elements*, the emblematic Proposition for the hypothesis of Geometric Algebra.

*Pythagorean theorem. II.4/5.*

If an orthogonal triangle has sides a, b around the right angle and hypotenuse c,
then $a^2 + b^2 = c^2$.
[Proof from Proposition II.4.]

Propositions II.5 & 6 are each equivalent to Proposition II.4, given in a completion of the square form. With the Pythagorean theorem and the completion of squares II.5 & II.6 available, it is quite simple to prove the two fundamental forms of the Pythagorean Application of Areas:

*Proposition II.5.* If a, x are lines, with $a/2 > x$, then $(a/2)^2 - (a/2 - x)^2 = x \cdot (a - x)$.

*Pythagorean Application of Areas in Defect. II.5/6.*
Given lines a, m, with $a/2 > m$,
to construct a line x, such that $x \cdot (a - x) = m^2$.
[Proof from Pythagorean theorem and Proposition II.5.]

*Proposition II.6.* If a, x are lines, then $(a/2 + x)^2 = (a/2)^2 + x \cdot (a + x)$.

*Pythagorean Application of Areas in Excess. II.6/7.*
Given lines a, m, to construct a line x, such that $x \cdot (a + x) = m^2$.
[Proof from Pythagorean theorem and Proposition II.6: construct the triangle with sides m, a/2 about a right angle and set y for the hypotenuse of the triangle.
By the Pythagorean theorem $y^2 = m^2 + (a/2)^2$. Certainly $y > a/2$. Set $x = y - (a/2)$.
Apply II.6: $y^2 = (a/2 + x)^2 = (a/2)^2 + x \cdot (a + x)$. Hence $x \cdot (a + x) = m^2$.]

The close relation of the construction of a mean proportional of two lines [II.14] with the Pythagorean Application of Areas in Defect [II.5/6] is seen by the following

*Proposition.* Let a be a given line.
(a) [II.5/6] if a line m is given, such that $a/2 > m$, then construct line x, such that $x \cdot (a - x) = m^2$;
(b) [II.14] if a line x is given, such that $a/2 > x$, then construct line m, such that $x \cdot (a - x) = m^2$.

*Proof.* By the Pythagorean theorem II.4/5 and Proposition II.5, applied in that order for (a) and in reverse order for (b).

**6.2.** *The use of the Pythagorean Application of Areas for the generation of the infinite anthyphairesis*

**6.2.1.** *For the computation of*
*the infinite anthyphairesis of a to b, satisfying an equation of the form*
$A \cdot a^2 = C \cdot b^2$,
*it is necessary to deal with*
*the computation of the anthyphairesis of a to b satisfying an equation of the form*
$A_1 \cdot b^2 = B_1 \cdot b \cdot c + C_1 \cdot c^2$

The computation of the anthyphairesis of a to b, in case $A \cdot a^2 = C \cdot b^2$, with A, C natural numbers reduces to the computation of the anthyphairesis of b to c, $\text{Anth}(a, b) = [k, \text{Anth}(b, c)]$, satisfying a quadratic equation in excess $A_1 \cdot b^2 = B_1 \cdot b \cdot c + C_1 \cdot c^2$.

It is easy to see that the proof of the anthyphairesis of a, b, a > b, with $A \cdot a^2 = C \cdot b^2$ for non-square natural numbers A, C, leads necessarily, from the second step on, to the proof of the incommensurability of the quadratic equations in excess of the form $A_1 \cdot a^2 = B_1 \cdot a \cdot b + C_1 \cdot b^2$, for natural numbers $A_1$, $B_1$, $C_1$, with non-square discriminant $A_1^2 + 4 \cdot B_1 \cdot C_1 = 4 \cdot A \cdot C$.

*Proposition.*
If a, b are lines, with a > b, and
$A \cdot a^2 = C \cdot b^2$, with A, C natural numbers with A < C, and non-square $A \cdot C$
[in particular if
$a^2 = N \cdot b^2$ for some non-square natural number N],
then there are a number k, a line c, c < b, and natural numbers $A_1$, $B_1$, $C_1$, such that
$\text{Anth}(a, b) = [k, \text{Anth}(b, c)]$, and $A_1 \cdot b^2 = B_1 \cdot b \cdot c + C_1 \cdot c_1^2$.

**Proof.** Let J be a natural number such that $J^2 < 4 \cdot A \cdot C < (J + 1)^2$, and k the integer part of $J/2 \cdot A$. Then $a = k \cdot b + c$, c < b, $2 \cdot A \cdot k < J$, $A \cdot k^2 < C$.
Set $A_1 = C - A \cdot k^2$, $B_1 = 2 \cdot A \cdot k$, $C_1 = A$, and verify that
$A_1 \cdot b^2 = B_1 \cdot b \cdot c + C_1 \cdot c^2$, and $\text{Anth}(a, b) = [k, \text{Anth}(b, c)]$.

*Thus, even if we are interested mainly or only for the incommensurability of the case $a^2 = N \cdot b^2$, we are forced to consider the more general quadratic equation $A \cdot a^2 = B \cdot a \cdot b + C \cdot b^2$.*

**6.2.2**. *The anthyphairetic division $a = k \cdot b + c$ of line a to line b, satisfying a quadratic equation $A \cdot a^2 = B \cdot a \cdot b + C \cdot b^2$ is computed by means of the Pythagorean method of Application of Areas; and the anthyphairetic substitution results in a quadratic equation $A_1 \cdot b^2 = B_1 \cdot b \cdot c + C_1 \cdot c^2$*

*Proposition.*
(a) An arithmetical quadratic equation in excess
$\quad A \cdot a^2 = B \cdot a \cdot b + C \cdot b^2$,
where a, b are lines, with a > b and A, B, C are natural numbers with non-square discriminant $B^2 + 4 \cdot A \cdot C$, is transformed,
upon setting $x = A \cdot a - B \cdot b$,
into the *Application of Areas in excess*
$\quad x \cdot (B \cdot b + x) = A \cdot C \cdot b^2$.
(b) The solution of the Application of Areas in excess yields, upon setting
$\quad$ J the natural number such that $J^2 < B^2 + 4 \cdot A \cdot C < (J + 1)^2$, and
$\quad$ k the integer part of $(B + J) / 2 \cdot A$,
the anthyphairetic relation
$\quad a = k \cdot b + c$, with c < b.
(c) $B < 2 \cdot A \cdot k \leq J + B$, $A \cdot k^2 < B \cdot k + C$.
(d) Upon setting
$\quad A_1 = B \cdot k + C - A \cdot k^2$, $B_1 = 2 \cdot A \cdot k - B$, $C_1 = A$,
we obtain again an arithmetical quadratic equation in excess
$\quad A_1 \cdot b^2 = B_1 \cdot b \cdot c + C_1 \cdot c^2$.

*Proof.* (a) is routine.
(b) The solution of the Application of Areas, according to Proposition II.6/7, yields
[b1] $2 \cdot y = B \cdot b + 2 \cdot x$
[b2] $2 \cdot y = J \cdot b + c_1$, with $c_1 < b$.
From [b1] and [b2] we get
[b3] $2 \cdot x = (J - B) \cdot b + c_1$, with $c_1 < b$.
Since $x = A \cdot a - B \cdot b$, we get
[b4] $2 \cdot A \cdot a = 2 \cdot B \cdot b + 2 \cdot x = (B + J) \cdot b + c_1$, with $c_1 < b$.
[b5] $a = ((B + J)/2 \cdot A) \cdot b + c_1 / 2 \cdot A$
[b6] $a = k + (D \cdot b + c_1)/2A$, with $D < 2 \cdot A$.
Setting $(D \cdot b + c_1)/2 \cdot A = c$, and noting that $c < b$, we conclude that
[b7] $a = k \cdot b + c$, with $c < b$.
Thus k is the integer part of the arithmetical ratio $(J + B)/2 \cdot A$
[This is the general algebraic step, obtained from
the solution of the Application of Areas in Excess]

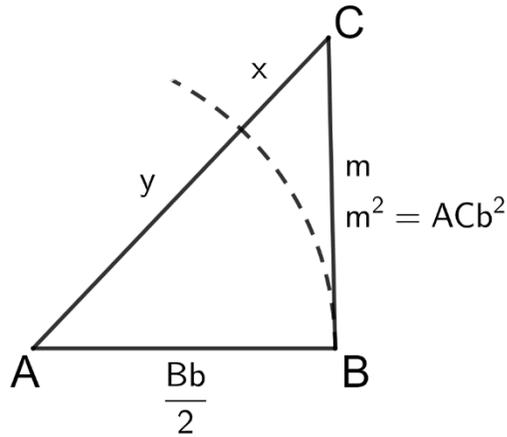

(c) Clearly $a < (k + 1) \cdot b \leq 2 \cdot k \cdot b$ and $B \cdot b < A \cdot a$, hence $B \cdot b < A \cdot a < 2 \cdot A \cdot k \cdot b$,
hence $B < 2 \cdot A \cdot k$.
By the definition of k, $2 \cdot A \cdot k \leq J + B$, hence $2 \cdot A \cdot k - B \leq J$.
By the definition of J, we have $(2 \cdot A \cdot k - B)^2 \leq J^2 < B^2 + 4 \cdot A \cdot C$,
hence by II.7, we have $4 \cdot A^2 \cdot k^2 + B^2 < B^2 + 4 \cdot A \cdot C + 4 \cdot A \cdot B \cdot k$, $4 \cdot A^2 \cdot k^2 < 4 \cdot A \cdot C + 4 \cdot A \cdot B \cdot k$,
hence $A \cdot k^2 < B \cdot k + C$.
(d) By substitution $A \cdot (k \cdot b + c)^2 = B \cdot (k \cdot b + c) \cdot b + C \cdot b^2$,
hence $A \cdot k^2 \cdot b^2 + A \cdot c^2 + 2 \cdot A \cdot k \cdot b \cdot c = B \cdot k \cdot b^2 + B \cdot b \cdot c + C \cdot b^2$,
hence $(B \cdot k + C - A \cdot k^2) \cdot b^2 = (2 \cdot A \cdot k - B) \cdot b \cdot c + A \cdot c^2$,
hence $A_1 \cdot b^2 = B_1 \cdot b \cdot c + C_1 \cdot c^2$.

It is next proved that the quadratic equation in defect $A \cdot a^2 + C \cdot b^2 = B \cdot a \cdot b$ is necessarily transformed into a quadratic equation in excess, after a finite number of successive steps of anthyphairetic substitution.

**6.2.3.** *Proposition.* An arithmetical quadratic equation
either in defect $A \cdot a^2 + C \cdot b^2 = B \cdot a \cdot b$,

or of the form $A \cdot a^2 + B \cdot a \cdot b = C \cdot b^2$,
or of the form $A \cdot a^2 = C \cdot b^2$
with lines a, b incommensurable,
under repeated anthyphairetic substitution eventually turns into an arithmetical quadratic equation in excess.

*Proof.*
Let $A \cdot a^2 + C \cdot b^2 = B \cdot a \cdot b$ be a quadratic equation in defect, with a, b incommensurable, $a > b$.
Then the discriminant $B^2 - 4 \cdot A \cdot C$ is a non-square number. We set $x = B \cdot b - A \cdot a$, and obtain the Application of Areas in defect II.5/6, $x \cdot (B \cdot b - x) = A \cdot C \cdot b^2$.
We set J the natural number such that $J^2 < B^2 - 4 \cdot A \cdot C < (J + 1)^2$,
and k the integer part of $(B + J) / 2A$.
The solution of the Application of Areas in defect yields, in analogy to the corresponding step for a quadratic equation in excess (Proposition 5.2.2), the anthyphairetic relation $a = k \cdot b + c$, with $c < b$.
[This is the general algebraic step]

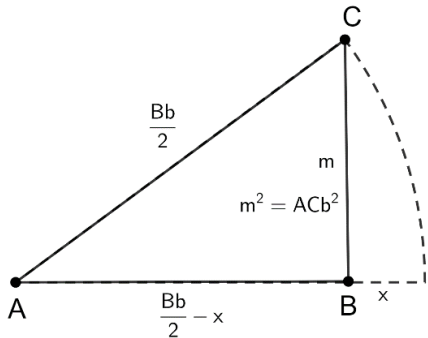

Next we employ the anthyphairetic substitution of a by $k \cdot b + c$ in the initial equation $A \cdot a^2 + C \cdot b^2 = B \cdot a \cdot b$, and we have $A \cdot (k \cdot b + c)^2 + C \cdot b^2 = B \cdot (k \cdot b + c) \cdot b$;
by Propositions II.1, 4 we have
$(A \cdot k^2 + C) \cdot b^2 + 2 \cdot A \cdot k \cdot b \cdot c + A \cdot c^2 = B \cdot k \cdot b^2 + B \cdot b \cdot c$.

*Claim.* $A \cdot k^2 + C \neq B \cdot k$
[If $A \cdot k^2 + C = B \cdot k$, then $2 \cdot A \cdot k \cdot b \cdot c + A \cdot c^2 = B \cdot b \cdot c$, then $2 \cdot A \cdot k \cdot b + A \cdot c = B \cdot b$,
hence $A \cdot c = (B - 2 \cdot A \cdot k) \cdot b$, hence b, c commensurable, hence a, b commensurable, impossible]
We then distinguish the following cases:

| Case 1. | Case 2. | Case 3. | Case 4. | Case 5. |
|---|---|---|---|---|
| $B \cdot k > A \cdot k^2 + C$, $2 \cdot A \cdot k > B$ | $B \cdot k > A \cdot k^2 + C$, $2 \cdot A \cdot k < B$ | $B \cdot k > A \cdot k^2 + C$, $2 \cdot A \cdot k = B$ | $B \cdot k < A \cdot k^2 + C$, $2 \cdot A \cdot k < B$ | $B \cdot k < A \cdot k^2 + C$, $2 \cdot A \cdot k \geq B$ |
| $A_1 = B \cdot k - A \cdot k^2 - C$, $B_1 = 2 \cdot A \cdot k - B$, $C_1 = A$ are natural numbers | $A_1 = B \cdot k - A \cdot k^2 - C$, $B_1 = B - 2 \cdot A \cdot k$, $C_1 = A$ are natural numbers | $A_1 = B \cdot k - A \cdot k^2 - C$, $[B_1 = 0]$ $C_1 = A$ are natural numbers | $A_1 = A \cdot k^2 + C - B \cdot k$, $B_1 = B - 2 \cdot A \cdot k$, $C_1 = A$ are natural numbers | $A_1 = A \cdot k^2 + C - B \cdot k$, $B_1 = 2 \cdot A \cdot k - B$, $C_1 = A$ are natural numbers |
| The anthyphairetic substitution $a = k \cdot b + c$ in $A \cdot a^2 + C \cdot b^2 = B \cdot a \cdot b$ yields: | | | | |
|  | $A_1 \cdot b^2 + B_1 \cdot b \cdot c = C_1 \cdot c^2$ | $A_1 \cdot b^2 = C_1 \cdot c^2$ | $A_1 \cdot b^2 + C_1 \cdot c^2 = B_1 \cdot b \cdot c$ | Impossible |
|  |  |  | Application of Areas in defect with $B_1 < B$ |  |
|  | The next step: $b = m \cdot c + d, d < c$ $A_2 = C_1 - A_1 \cdot m^2 - B_1 \cdot m$ $B_2 = B_1 + 2 \cdot A_1 \cdot m$ $C_2 = A_1$ | The next step: $b = m \cdot c + d, d < c$ $A_2 = C_1 - A_1 \cdot m^2$ $B_2 = 2 \cdot A_1 \cdot m$ $C_2 = A_1$ |  |  |
| $A_1 \cdot b^2 = B_1 \cdot b \cdot c + C_1 \cdot c^2$ Application of Areas **in Excess** | $A_2 \cdot c^2 = B_2 \cdot c \cdot d + C_2 \cdot d^2$ Application of Areas **in Excess** | $A_2 \cdot c^2 = B_2 \cdot c \cdot d + C_2 \cdot d^2$ Application of Areas **in Excess** | Eventually Application of Areas **in Excess** |  |

The (many but manageable) details are left to the reader.

**Section 7.** *The anthyphairetic reconstruction of Theaetetus' proof of quadratic incommensurabilities (*Theaetetus, *147d3-e1), Part II: Theaetetus' New Idea, a Pigeonhole Argument, results in the preservation of the Application of Areas/Gnomons, hence both in Anthyphairetic Periodicity and in the Birth of the theory of proportion of magnitudes*

We present the basic new idea of Theaetetus, a *pigeonhole argument* based on the *invariance of the discriminant under anthyphairetic substitution* (Section 6.2). The argument implies that there are two stages in the anthyphairesis of an arithmetical quadratic equation, say the mth and afterwards the nth stage, so that the mth, say $e_m$, and the (m+1)th, $e_{m+1}$, anthyphairetic remainder satisfy the same arithmetical quadratic equation, satisfied by the nth $e_n$ and (n+1)th $e_{n+1}$ remainders, say

$A \cdot e_m^2 = B \cdot e_m \cdot e_{m+1} + C \cdot e_m^2$ and $A \cdot e_n^2 = B \cdot e_n \cdot e_{n+1} + C e_n^2$ (for some natural numbers A, B, C).
This is a generalization of the Pythagorean *Invariance of Gnomons*

As it is not difficult to see (shown in section 6.1, below) this repetition implies that

$$\text{Anth}(e_m, e_{m+1}) = \text{Anth}(e_n, e_{n+1}).$$

It is precisely this condition that, by some natural inductive argument,
the anthyphairesis of the initial line segments a to b is eventually periodic.

At this point, we compare

> Plato's imitation, involving Logos, in effect the Logos Criterion for anthyphairetic periodicity, e.g. B6/A6 = B9/A9 for the Angler (in Section 4) and
> the condition on equal anthyphairesis, which arose naturally from our reconstruction.

By this comparison we conclude that

> Theaetetus reached at the definition of proportion of magnitudes in terms of equal anthyphairesis, by interpreting the condition $\text{Anth}(e_m, e_{m+1}) = \text{Anth}(e_n, e_{n+1})$ within the anthyphairesis to denote proportion $e_m/e_{m+1} = e_n/e_{n+1}$.
> This marks
> *the birth of the Logos Criterion for anthyphairetic periodicity* (Section 6.3),
> Plato's dialectical meaning of his term *Logos*
> (in the *Sophist* as in Section 4 and in his whole work), and
> *the birth of Theaetetus' theory of ratios of magnitudes*.

Thus, Theaetetus' theory of ratios of magnitudes is, from its initial conception, inseparably connected with periodic anthyphairesis.

## 7.1. *The transition from the preservation of Gnomons to the equality of anthyphairesis*

The Propositions of this Section depend on

### 7.1.1. *Proposition I.44.*
Given lines a, b, c, to construct a line d, such that $a \cdot b = c \cdot d$.

We note that Proposition I.44 is in the second part of Book I of the *Elements*, which may be considered preparatory for Book II; and Book II culminates with the proof of incommensurability at Proposition II.8/9.

### 7.1.2. *Lemma.* If $A \cdot a^2 = B \cdot a \cdot b + C \cdot b^2$ and $A \cdot a^2 = B \cdot a \cdot d + C \cdot d^2$, then $b = d$.

**Proof.** By hypothesis $B \cdot a \cdot b + C \cdot b^2 = B \cdot a \cdot d + C \cdot d^2$. Suppose b is different from d, say $b > d$. Then $B \cdot a \cdot b > B \cdot a \cdot d$ and $C \cdot b^2 > C \cdot d^2$, hence $B \cdot a \cdot b + C \cdot b^2 > B \cdot a \cdot d + C \cdot d^2$, a contradiction.

### 7.1.3. *Proposition.* If $A \cdot a^2 = B \cdot a \cdot b + C \cdot b^2$ and $A \cdot c^2 = B \cdot c \cdot d + C \cdot d^2$, then $a \cdot d = b \cdot c$.
The condition of the Proposition is referred to as
*the preservation of quadratic equations/Gnomons.*

**Proof.** Let r be an assumed line.

By Proposition I.44, there is a line s, such that a·r = c·s.
By Proposition I.44, there is a line e, such that b·r = e·s.
From the assumption we have A·a$^2$·r = B·a·b·r + C·b$^2$·r,
whence A·a·c·s = B·a·e·s + C·b·e·s,
whence A·a·c = B·a·e + C·b·e,
whence A·a·c·r = B·a·e·r + C·b·e·r,
whence A·c$^2$·s = B·c·e·s + C·e$^2$·s,
whence A·c$^2$ = B·c·e + C·e$^2$.
By Lemma 7.1.2, d = e. Hence, a·r = c·s and b·r = d·s.
But then (a·r)·d = (c·s)·d = c·(s·d) = c·(b·r), hence (a·d)·r = (c·b)·r, hence a·d = c·b.

**7.1.4.** *Proposition.* If A·a$^2$ = B·a·b + C·b$^2$ and A·c$^2$ = B·c·d + C·d$^2$, then Anth(a, b) = Anth(c, d).

*Proof.* If A·a$^2$ = B·a·b + C·b$^2$, and A·c$^2$ = B·c·d + C·d$^2$, then a·d = b·c (Proposition 7.1.3); and, if a·d = c·b, then Anth(a, b) = Anth(c, d) (Proposition 5.2.1).

**7.2.** *Theaetetus' fundamental new idea: the pigeon-hole argument.*
*According to our interpretation, the mathematical idea that Theaetetus had, according to Plato's* Theaetetus *147d7, was a pigeonhole principle, based on the invariance of the discriminant of the quadratic equation under anthyphairetic substitution*

**7.2.1**. Propositions 7.1.3 and 7.1.4 are natural extensions of the Pythagorean and Theodorus' approach and they do not qualify as the major new idea, described in the *Theaetetus* 147d7. The cases worked out by Theodorus, a$^2$ = N·b$^2$ for non-square N up to N = 17, lead to the formation of the discriminant and do indicate the reappearance of the same discriminant in every step of the solution of the required Application of Areas in Excess, suggesting the invariance of the discriminant under anthyphairetic substitution.
The fundamental new idea that Theaetetus must have had as a result of Theodorus' lesson, can only be the *pigeon-hole principle* that is needed for the proof of the general periodicity theorem. Theaetetus must have noticed that in all the cases of the incommensurability proofs by Theodorus the discriminant of the quadratic equation in excess remains constant.
As noted in detail in Negrepontis, Farmaki, Brokou, 2024, Section 6.5, Plato's appeal, in the *Politicus* 272d-e, for the necessity of transferring from the Kronos era to the Zeus era in his related myth, appears to be a philosophical version of a pigeonhole argument, employed exactly at the point where it is needed in the proof of Theaetetus' palindromic periodicity theorem, a theorem deeper than the simple periodicity theorem presently considered and which must have been proved at a later stage.

**7.2.2**. *The preservation of the discriminant of an arithmetic quadratic equation under anthyphairetic substitution*

*Proposition.* If lines a, b, with a > b, satisfy the quadratic equation in excess A·a$^2$ = B·a·b + C·b$^2$, where A, B, C are natural numbers with non-square discriminant B$^2$ + 4·A·C,
a = k·b + c, with c < b, and A$_1$ = B·k + C – A·k$^2$, B$_1$ = 2·A·k – B, C$_1$ = A, as in Proposition 6.2.2,

then the discriminant of the produced quadratic equation in excess $A_1 \cdot b^2 = B_1 \cdot b \cdot c + C_1 \cdot c^2$
is equal to the discriminant of the given quadratic equation in excess.

The *Proof* is routine:
$B_1^2 + 4 \cdot A_1 \cdot C_1 =$
$= (2 \cdot A \cdot k - B)^2 + 4 \cdot (B \cdot k + C - A \cdot k^2) \cdot A =$ (by Proposition II.7)
$(4 \cdot A^2 \cdot k^2 + B^2 - 4 \cdot A \cdot B \cdot k) + (4 \cdot A \cdot B \cdot k + 4 \cdot A \cdot C - 4 \cdot A^2 \cdot k^2) =$
$B^2 + 4 \cdot A \cdot C$.

**7.2.3.** *Theaetetus' general periodicity theorem*. If the lines a, b, a > b, satisfy
**either** $A \cdot a^2 = C \cdot b^2$ with non-square discriminant $4 \cdot A \cdot C$,
**or** a quadratic equation **in Excess** $A \cdot a^2 = B \cdot a \cdot b + C \cdot b^2$ with non-square discriminant $B^2 + 4 \cdot A \cdot C$
(including the case $A \cdot a^2 = C \cdot b^2$)
**or** a quadratic equation **in Defect** $A \cdot a^2 + C \cdot b^2 = B \cdot a \cdot b$ with non-square discriminant $B^2 - 4 \cdot A \cdot C$,
then the anthyphairesis of a to b is **eventually periodic**.

*Proof.* We denote the anthyphairesis of a to b as follows:
$a = k_0 \cdot b + e_1$, with $b > e_1$,
$b = k_1 \cdot e_1 + e_2$, with $e_1 > e_2$,
…
$e_{m-1} = k_m \cdot e_m + e_{m+1}$, with $e_m > e_{m+1}$,
…
$e_{n-1} = k_n \cdot e_n + e_{n+1}$, with $e_n > e_{n+1}$,
…
A quadratic equation in defect turns into a quadratic equation in Excess after a finite number of anthyphairetic substitutions, according to Proposition 6.2.3. The discriminant is kept invariant in all anthyphairetic substitutions by Proposition 7.2.2.
*A pigeonhole argument,* based on the invariance of the discriminant, yields two equations in excess in steps, say m < n, with the same integer coefficients:
$A_m \cdot e_m^2 = B_m \cdot e_m \cdot e_{m+1} + C_m \cdot e_{m+1}^2$ and
$A_n \cdot e_n^2 = B_n \cdot e_n \cdot e_{n+1} + C_n \cdot e_{n+1}^2$,
with $A_m = A_n$, $B_m = B_n$, $C_m = C_n$.
By Proposition 7.1.3, $e_m \cdot e_{n+1} = e_{m+1} \cdot e_n$.
By Proposition 5.2.1, $\text{Anth}(e_m, e_{m+1}) = \text{Anth}(e_n, e_{n+1})$.
It now follows by the Theaetetean definition of proportion that $e_m/e_{m+1} = e_n/e_{n+1}$ and therefore, by proposition 2.3.1, the anthyphairesis of a to b is eventually periodic.

*Note*. We believe that the birth of Theaetetus' theory of ratios of magnitudes occurred precisely at the moment when Theaetetus rephrased the established condition of the Pythagoreans regarding the preservation of Application of Areas/Gnomons, implying anthyphairetic periodicity and appearing at the conclusion of his theorem, in terms of equality of anthyphaireses and then in terms of proportion, equality of ratios, according to the definition provided to us by Aristotle in Topics 158b29-35. Thus anthyphairetic periodicity is built in the initial conception of the theory.

**8**. *Theaetetus' theory is **only** about ratios of magnitudes with eventually periodic anthyphairesis*

**8.1**. *The birth of Theaetetus' theory of ratios of magnitudes*

According to our reconstruction of Theaetetus' theorem, reported in Plato's Theaetetus 147d-148b, in Section 7, the birth of Theaetetus' theory of ratios of magnitudes occurred at the exact point where the older, originally Pythagorean finitization of the infinite anthyphairesis by the preservation of Gnomons/Application of Areas was replaced by the definition of proportion in terms of the equality of anthyphairesis and the resulting Logos Criterion for the periodicity and the finitization of the knowledge of a geometric anthyphairesis. Thus the initial role of Theaetetus' theory of ratios was to rephrase the Pythagorean finitizer, the preservation of the Gnomons, which strictly speaking applied only for the infinite anthyphairesis of the diameter to the side of a square with a more general finitizer, the Logos Criterion, that applied to very eventually periodic anthyphairesis.

**8.2**. *Logos in Plato's philosophy*

According to our interpretation of Plato's Logos, as used in the Sophist and in the Meno, in Section 5, but in fact throughout Plato's philosophy, Plato's Logos is the philosophical imitation of the Logos Criterion, the fundamental tool of acquiring knowledge of a Platonic intelligible Being, Plato's Idea, which intelligible Being therefore has the structure of a dyad of kinds in imitation of a dyad of lines in geometric periodic anthyphairesis. Thus Platonic Logos is the philosophic finitizer and provider of Knowledge of the infinite anthyphairetic structure of a Platonic Idea.
Thus Theaetetus' theory of ratios of magnitudes was conceived as a tool to study eventually periodic anthyphairesis, and Plato used Theaetetus' theory in order to describe his philosophical theory of intelligible Beings. Theaetetus' theory of ratios turned out to be crucial for Plato's philosophical system, because this theory was initially conceived by the Eleatics Parmenides and Zeno who had no theory of ratios of magnitudes and no Logos to rely upon. Zeno's finitizer of the philosophical infinite anthyphairesis, described in Fragments B3 & B1, and mentioned in the Platonic dialogues Parmenides Second Hypothesis and Sophist, is the precursor of Platonic Logos (as shown in Negrepontis, TA).
If the basic role of Theaetetus' theory of ratios is to be an effective finitizer and producer of knowledge and Oneness, by applying it to eventually periodic anthyphairesis, equivalently to quadratic equations in the form of Pythagorean Application of Areas with arithmetical coefficients, then Theaetetus' theory cannot be a theory for all ratios of magnitudes, but only for ratios of magnitudes a and b, such that the anthyphairesis of a to b is eventually periodic.

**8.3**. *Definition 3 in the Theaetetean Book X of the* Elements

That Theaetetus' theory is a theory of ratios only for ratios in eventually periodic anthyphairesis may be seen from the Definition X.3 of the *Elements* on rational and alogos/irrational lines, and from the remarkable *Anonymon Scholion X.2 to the Elements.*
Theaetetus' theory of ratios is applied to Book X of the *Elements*, which is generally considered to be due to Theaetetus, and concerns the theory of "alogoi"/irrational lines.

Definition X.3 of the *Elements*.
Let then the assigned straight line be called rational (rhete), and
those straight lines which are commensurable with it,
whether in length and in square, or in square only, rational (rhetai),
but those that are incommensurable with it irrational (alogoi).

Thus, given an assumed line r, a line a possesses ratio according to the Definition X.3 with respect to r if a and r are commensurable in square, namely if $a^2$, $r^2$ are commensurable squares. But according to Theaetetus' theorem, the anthyphairesis of a to r is eventually periodic. Thus the condition that a line a possesses ratio with respect to a line r in Definition X.3 clearly refers to the Theaetetus' theory of ratios, and is a condition closely related and in fact stronger than eventual anthyphairesis. (Two *Anonyma Scholia to the Elements*, *Scholion* V.22 and *Scholion* V.24 refer to this definition).

**8.4**. *The comparison of Theaetetus' to Eudoxus' theory of ratios of magnitudes in the* Anonymon Scholion X.2 to Euclid's Elements *and Proclus,* Commentary to Euclid *285,5-286,11*

The situation is quite different for Eudoxus' theory of ratios of magnitudes, presented in Book V of the *Elements*.

Definition V. 4 of the *Elements*
Magnitudes are said to have a ratio to one another
which can, when multiplied, exceed one another.

Thus Eudoxus' ratio of magnitudes is defined for those magnitudea a and b, such that Eudoxus' condition (mentioned also in Section 3) is satisfied, namely a>nb and b>ma for some numbers m and n.

The remarkable *Scholion X.2 to the Elements* compares the two definitions, the Theaetetean and the Eudoxean and criticizes the Eudoxean theory,
presented in Book V of the *Elements*,
as being so general that it annihilates the difference between rational and alogoi lines,
defined in the Theaetetean Book X of the *Elements*.
[we omit the references, in *Scholion X.2*, to commensurable lines, to make clearer the opposition between the two theories]

*Scholion X.2 to the Elements*

| |
|---|
| Τὰ […] μαθήματα *φανταστικῶς* νοοῦμεν, |
| We obtain knowledge of the magnitudes by imagination. |
| […] |
| διὸ καὶ τὰ [μαθήματα] εἰς ἄπειρον διαιρεῖται, |
| therefore the magnitudes are divided *ad infinitum* |
| […] |
| ἡ δὲ *φαντασία* πλῆθος ἄπειρον ἔχει· διὸ τὰ *φανταστὰ* ἄπειρα. |

> Imagination possesses an infinite multitude; thus entities of imagination are infinite.
> καὶ τὰ μεγέθη οὖν ὡς φανταστὰ ἄπειρα καὶ ἡ τομὴ αὐτῶν.
> And the magnitudes as being entities of imagination are infinite and infinitely divisible
> *Scholion X.2, lines 1-6*

The unusual term "phantasia" refers to the divisibility *ad infinitum* of the magnitudes. We conjecture that the origin of the term is precisely in the sentence 147d7-8 of the *Theaetetus*, "since the powers "ephainonto" infinite in multitude", and was probably concocted by Proclus. There is a long passage in Proclus' *Commentary to Euclid*, 285,5-286,11 about imagination and the infinite (cf. also *Commentaries to Euclid* 50-58, where imagination is also called "pathetikos nous", passive intellect):

> *ἐν τῇ φαντασίᾳ τὸ ἄπειρον ὑφίστασθαι μόνον*
> οὐ νοούσης *τὸ ἄπειρον* τῆς φαντασίας.
> in the imagination the *infinite* only exists
> without the imagination's knowing *the infinite* 285,5-7
> ὥσπερ γὰρ τὸ σκότος τῷ μὴ ὁρᾶν ἡ ὄψις γινώσκει,
> οὕτως ἡ φαντασία *τὸ ἄπειρον* τῷ μὴ νοεῖν.
> Just as the sight recognizes the darkness
> by not seeing,
> so the imagination (φαντασία) recognizes the *infinite* (ἄπειρον)
> by not understanding it.

> ἅμα γὰρ νοεῖ
> καὶ μορφὴν ἐπάγει τῷ νοουμένῳ καὶ πέρας,
> καὶ τῇ νοήσει
> τὴν τοῦ φαντάσματος ἵστησι διέξοδον,
> καὶ διέξεισιν αὐτὸ
> καὶ περιλαμβάνει. 285,7-10
> For when the imagination knows,
> it simultaneously assigns to the object of its knowledge a Form and a Finite,
> and in knowing
> brings to an end the movement through the imagined object, and
> goes through it, and
> comprehends it.

But when the finitizer (peras), namely the Logos, is added, then
the knowledge of the infinite is obtained, and
the infinite is comprehended (perilambanei) (into One).
We thus understand that the meaning of the *imagination* is
the *(anthyphairetic) infinite without Logos and knowledge.*

*Scholion X.2* in the sequel compares Theaetetus' and Eudoxus' theory of ratios:

| *possessing Logos* | *possessing Logos* |
|---|---|

| *in the sense of Eudoxus' theory* | *in the sense of Theaetetus' theory* |
|---|---|
| Εἰ πάντα τὰ μεγέθη τὰ πεπερασμένα δύναται πολλαπλασιαζόμενα ἀλλήλων ὑπερέχειν (τοῦτο δὲ ἦν **τὸ λόγον ἔχειν,** ὡς ἐν τῷ πέμπτῳ μεμαθήκαμεν), if all finite magnitudes have the power by being multiplied to exceed each other (for such was the condition that a magnitude *possesses logos,* as we were taught in Book Five of the *Elements*) | |
| *τίς μηχανὴ τὴν τῶν ἀλόγων ἐπεισφέρειν διαφοράν; how on earth (tis mechane) [Eudoxus theory] can ever distinguish between the rational and the alogon?* | |
| | ἢ ὅτι τὸ μέτρον […] θέσει δὲ ἐν τοῖς μεγέθεσι διὰ τὴν ἐπ' ἄπειρον τομήν; *πρὸς γὰρ πῆχυν ἢ σπιθαμὴν* ἤ τι τοιοῦτον *γνώριμον* μέτρον τὸ ῥητὸν καὶ τὸ ἄρρητον *γνωρίζομεν.* while we obtain knowledge of the rational and the irrational with respect to a known assumed line [*definition X.3, Theaetetus theory*] |
| καὶ μὴν **τὸ λόγον ἔχειν** to possesse logos is meant ||
| ἄλλως μὲν ἐπὶ τῶν μεγεθῶν λέγεται τῶν πεπερασμένων καὶ ὁμογενῶν, differently for magnitudes finite and homogeneous | |
| | […] ἄλλως ἐπὶ τῶν ῥητῶν προσαγορευομένων. differently for rational magnitudes |
| ὅπου μὲν γὰρ *ὁ λόγος μόνον* καὶ ἡ σχέσις θεωρεῖται τῶν πεπερασμένων μεγεθῶν κατὰ τὸ μεῖζον καὶ ἔλαττον, for [in Eudoxus theory] | |
| *logos* and the relation of the finite magnitudes is considered *only* according to the greater and the smaller [*definition V.4 of the* Elements] | |
| | […] ὅπου δὲ πρὸς τὸ ἐγκείμενον μέτρον τὴν τῶν ῥητῶν ἡμῖν πρὸς τὰ ἄλογα διαφορὰν παρέσχετο. while [in definition X.3 of the *Theaetetus' theory*] the assumed measure *does provide the difference between the ratiomal and the irrational lines* Scholion X.2, lines 7-22 |

.
The *Scholion* is strongly critical of Eudoxus' theory, exactly because it obliterates the difference between rational and irrational lines, literally between those lines that possess Logos and those that do not.
Thus Theaetetus' theory is clearly the superior one both for mathematical use and of course for Plato's philosophy, exactly because it is limited only to ratios that provide knowledge of the infinite, namely to ratios with eventually periodic anthyphairesis.

### 8.5. *Conclusion*

From 8.1-8.4, we have no doubt that, contrary to the reconstructions put forward by Becker (1933), van der Waerden (1954), and Knorr (1975) (Section 3), Theaetetus' theory of ratios of magnitudes is NOT about the class of all ratios of magnitudes, as Eudoxus' theory of ratios, presented in Book V of the Elements, BUT is strictly about the limited class of all ratios with finite or eventually periodic anthyphairesis.

### 9. *Our critique of the failure of earlier historians of Greek Mathematics and scholars of Plato to realize the fundamental role that the philosophical imitation of geometrical periodic anthyphairesis plays in Plato's philosophy, in particular in the dialogues* Theaetetus, Sophist *and* Meno*, and to reach an essential understanding of both Greek Mathematics and Plato's philosophy*

Bertand Russell, 1946, p.131 writes:

> It is noteworthy that modern Platonists, with few exceptions, are ignorant of mathematics in spite of the immense importance that Plato attached to arithmetic and geometry, and the immense influence that they had on his philosophy.

Paraphrasing Bertrand Russell, we may say that both modern Platonists and historians of Greek Mathematics did not realize the immense importance that Plato attached to periodic anthyphairesis, and the immense influence that it had on his philosophy. As a result Plato's philosophy was ill understood and some fundamental achievements of Greek Mathematics had lost the only source that could save them from oblivion.

Euclid with his proof of Proposition X.9 of the *Elements* may have contributed negatively towards a true understanding of the method that Theaetetus used for the proof of his general theorem, reported by Plato in the *Theaetetus* 147d-148b. Quite independently from the corrections that are needed to make the proof rigorous, our main concern is that Euclid gives an arithmetical proof, that is close to the methods used by Archytas, rather than those of Theaetetus (Section 9.1). Van der Waerden and Knorr might be misdirected by Euclid's proof of X.9 and *Scholion* X.62 towards the translations (Sections 9.2, 9.5) and the arithmetical proofs (Sections 9.3, 9.6) they have suggested.

Historians of Greek Mathematics and scholars on Plato did not take into account the revealing close connection between Theaetetus' proof of incommensurabilities, reported in the *Theaetetus* 147d-148b and the philosophical descriptions of the Knowledge of intelligible Beings, as imitations of Theaetetus' proof (148c-d), and as Name & Logos, equivalently as True Opinion of Logos (201-202), realized in the *Sophist* and in the *Meno*, along with other Platonic dialogues. Burnyeat cannot be right, either in the translation of the *Theaetetus* 147d-e (Sections 9.2, 9.5), or in the claim that Plato had no motive to indicate Theaetetus' method (Sections 9.4, 9.6). In our opinion the most serious failure of the Platonists is the misunderstanding of the anthyphairetic nature of Plato's Logos (Section 9.7).

**9.1.** *Theaetetus' theorem receives a defective but correctible proof in Proposition X.9 of the* Elements, *which is arithmetical, probably due to Archytas*

Theaetetus' theorem described in the *Theaetetus* 147d3-e1 passage has been given a problematic arithmetical proof rendered in different ways by various authors, and this difference has wide implications both for the history of Greek Mathematics and for Plato's philosophy. We will discuss the Euclidean defective arithmetical proof of the theorem of Theaetetus given in the *Elements* as Proposition X.9, provide a completion/correction of the proof, and indicate that it is an Archytan type proof. The *Anonymous Scholion to Euclid's Elements* X.62 correlates Proposition X.9 with the difficult *Theaetetus* 147d3-e1 passage stating:

> Τὸ θεώρημα τοῦτο Θεαιτήτειόν ἐστιν εὕρημα,
> καὶ μέμνηται αὐτοῦ ὁ Πλάτων ἐν Θεαιτήτῳ
> This theorem is a Theaetetean find,
> and Plato recalls it in the *Theaetetus*

*Proposition X.9 of the* Elements
> [1] The squares on straight lines commensurable in length
> have to one another the ratio which a square number has to a square number; and
> [2] squares which have to one another the ratio which a square number has to a square number also have their sides commensurable in length.
> [3] But the squares on straight lines incommensurable in length
> do not have to one another the ratio which a square number has to a square number; and
> [4] squares [of straight lines]
> which do not have to one another the ratio which a square number has to a square number also do not have their sides commensurable in length either.

We are only interested for [4] which may be stated as follows

> If a, b are lines, such that $a^2 / b^2 = C / A$ with AC non-square number, then a, b are incommensurable.

*Proof of [4] in the* Elements.
> Finally, let the square on *a* not have to the square on *b* the ratio which a square number has to a square number. I say that *a* is incommensurable in length with *b*.
> For, if *a* is commensurable with *b*, then [by [1]] the square $a^2$ on *a* has to the square $b^2$ on *b* the ratio which a square number $m^2$ has to a square number $n^2$.
> But it does not, therefore *a* is not commensurable in length with *b*.

The arithmetical proof is problematic in that it is incomplete on two counts.
Indeed, since we have assumed that a and b are commensurable,
it follows that there are numbers m, n and a line c, such that

$$a = m \cdot c, \quad b = n \cdot c.$$

Then obviously $a^2 = m^2 \cdot c^2$, $b^2 = n^2 \cdot c^2$.

*The first difficulty in the Euclidean proof.*
We certainly want to conclude that

$$m^2 \cdot c^2 / n^2 \cdot c^2 = m^2 / n^2,$$

a *mixed proportion* of the type $x / y = p / q$, where x, y are magnitudes and p, q are natural numbers. Euclid never supplies the required definition of a mixed proportion, but states and proves Propositions X.5 and 6, propositions involving a mixed proportion, as if he had such a definition. Thus Euclid's proofs are problematic. As van der Waerden, 1954, p.175, notes

> Indeed, Euclid's proofs contain logical errors for which we cannot hold Theaetetus responsible. Moreover these proofs depend on the theory of proportions developed in Book V, following Eudoxus, which was not yet known to Theaetetus.

On the other hand, it is quite natural to supply an adequate definition, most probably the one that Theaetetus himself had given:

> *Definition 9.1.1.* If a, b are magnitudes, and m, n are numbers, then the proportion $a / b = m / n$ holds if Anth(a, b) = Anth(m, n).

*Note*. Earlier definitions by Knorr (1975, pp. 253-254) and Taisbak (1982, pp. 21-25) provide an arithmetical definition of a mixed ratio and correct the proof of Proposition X.5.
With this definition it is "immediately obvious" (εὐθέως φανερὸν) as Aristotle would say (*Topics* 158b33), that $m^2 \cdot c^2 / n^2 \cdot c^2 = m^2 / n^2$, and hence that
$$a^2 / b^2 = m^2 / n^2.$$

Hence
$$m^2 / n^2 = C / A.$$
Note that this is a completely arithmetical relation

*The second difficulty in the Euclidean proof.*
The left-hand side is a ratio of the form square number to square number, while the right-hand side is not. Euclid appears to regard it as an obvious contradiction, but it certainly is not.
In order to show that this relation is impossible, we need to appeal to
Proposition VIII.8 of the *Elements*.

> If there are k mean proportionals between two numbers,
> then there are k mean proportionals between any two numbers
> which have the same ratios with the original numbers.

Indeed, since the terms of the ratio $m^2/n^2$ possess the mean proportional $m \cdot n$, and since $m^2 / n^2 = C / A$, it follows from Proposition VIII.8 that the numbers C, A possess also a mean proportional M, namely $AC = M^2$. But AC is not a square number, hence we now have a contradiction.

The problematic nature of the Euclidean proof of Proposition X.9 has been pointed out by van der Waerden, Knorr, and Mazur:

*van der Waerden, 1954*
> In Book X of the Elements, two proofs occur, so completely different in style from the rest of the book, that they have to be considered as later additions, viz. the proofs of X 9, 10. For, whereas Book X as a whole excels in its strictly logical structure and its extremely brief and elegant proofs, these two proofs not only contain superfluous parts, but also gaps. p. 159

*Knorr, 1975*
> Euclid does not offer a construction of roots in his demonstration of X,9, a theorem related to those in discussion in the *Theaetetus-passage.* p.71
>
> Theorem 7. Two lines are commensurable in length if and only if the squares constructed on them have the ratio which a square integer has to a square integer. *(Elements X,9)*], p.225
>
> Turning now to Theorem 7 we observe that although it has been stated in its Euclidean form (X,9),
> it is not in this form serviceable as a criterion of incommensurability.
> In addition, a result of the following type is needed:
> THEOREM 8. The least terms of a ratio of square integers are also square integers. p.232

*Mazur, 2007*
> Specifically, you will find—in various places—the claim that Theaetetus' theorem is proven in Proposition 9 of Book X of Euclid's *Elements*. It doesn't serve any purpose here to list the places where you find this incorrect assertion, except to say that it is incorrect, and it remains a thriving delusion since at least one important article published
as late as 2005 repeats it. It is an especially strange delusion since nothing subtle is going
on here.

Even a cursory glance at Proposition 9 of Book X will convince you that what is being
demonstrated there—if you take it in a modern perspective—is an utter triviality.
Proposition 9 of Book X stands, though, for an important issue in ancient thought if taken
on its own terms, but it won't prove irrationality of anything for us, let alone irrationality
of all the numbers that Theaetetus proves.

One might imagine that Heath's commentary on this
— which is perfectly clear, and says exactly what is indeed proved in Proposition 9 —
would dispel the misconception that Theaetetus' theorem about the irrationality of surds
is contained in this proposition, but it seems that this has held on with some tenacity.
p.236

The corrected proof is arithmetical, and thus not the original Theaetetus' proof, according to our argument in Sections 5-7, certainly not by Theaetetus. We note that Proposition VIII.8 of the *Elements* is attributed to Archytas, and thus the completed/corrected arithmetical proof of X.9 must be attributed to Archytas or his students. Since Book X is certainly attributed to Theaetetus, it follows that Euclid modified the original work.

## 9.2. *The translation of the Theaetetus 147d3-e1 passage by historians of Greek Mathematics van der Waerden and Knorr, and by Platonic scholar Burnyeat, 1978, and our critique of their translations.*

### 9.2.1. *Translation by van der Waerden, 1975, p.142 and p.166*

Here our Theodorus drew (or wrote) something about sides of squares
*(Περὶ δυνάμεών τι ἔγραφε)*
and showed *(ἀποφαίνων)* that those of three or five feet
are not commensurable in length with those of one foot,
and in this manner he took up one after another up to the one of 17 feet;
here something stopped him (or: here he stopped). p.142

Now it occurred to us (Theaetetus and another young man),
since the number of roots appeared to be infinite.
to try to collect them under one name,
by which we could henceforth call all the roots. p.166

### 9.2.2. *Translation by Knorr,1975, p.62-63*

Theaet.: Theodorus was proving for us via diagrams something about powers,
in particular about the three-foot-power and the five-foot-power –
demonstrating that these are not commensurable in length with the one-foot-power,
- and selecting each power individually in this way up to the seventeen-footpower;
but in this one for some reason he encountered difficulty.
Now, this is what occurred to us:

that, since we recognized the powers to be unlimited in number,
we might try to collect them under a single name,
by which we would designate all these powers.

**9.2.3.** *Translation by Burnyeat, 1978, p.493*

THEAETETUS. Theodorus here was drawing diagrams to show us something about
powers namely that a square of three square feet and one of five square feet
aren't commensurable, in respect of length of side, with a square of one square foot;
and so on, selecting each case individually, up to seventeen square feet.
At that point he somehow got tied up.
Well, since the powers seemed to be unlimited in number,
it occurred to us to do something on these lines:
to try to collect the powers under one term
by which we could refer to them all.

## 9.3. *The arithmetical reconstructions by van der Waerden and by Knorr of the proof of Theaetetus' theorem on quadratic incommensurabilities*

We will briefly describe the method that van der Waerden and Knorr suggested, in mutually substantial agreement, as the method that Theaetetus used for his proof of the quadratic incommensurabilities.

**9.3.1.** *Reconstruction by van der Waerden, 1954*

What is for us more important than these trivial definitions, is the proposition given at the
end, very briefly but nevertheless clearly:
*Line segments, which produce a square, whose area is an integer, but not a square
number, are incommensurable with the unit of length.*
...
How would the proof of such a proposition be made with classical methods?
Indirectly of course.
Supposing that two lines, which produce squares of areas n and 1, had commensurable
lengths, one would try to show that *n* is a square number.
The proof necessarily breaks up into a geometrical and a number-theoretical part.
From the assumed commensurability of the sides, one first deduces geometrically
that the ratio of the sides is equal to that of two numbers *p* and *q,* and
that the ratio of the areas of the squares which they produce is equal to that of $p^2$ and $q^2$,
so that
(1) $n:1 = p^2:q^2$,
or, equivalently,
(2) $p^2=nq^2$
Then it has to be shown, number-theoretically, that it follows from (1) and (2)
that n is a square number. p.166-167

It is Zeuthen's judgment [Zeuthen, 1910. p. 395]. that the greatest merit of Theaetetus is however
not to be found in the geometric analysis of the irrationality problem,
but in the arithmetic part, in the proof therefore that (2) can hold only if $n$ is a square.
According to this Judgment, Theaetetus had discovered and proved several propositions from the arithmetical Books VII and VIII which were needed for this proof.

I do not share Zeuthen's view on this. For, as we have seen before, Book VII
is of older date and forms the foundation of the Pythagorean theory of numbers.
But, once one has at his disposal the propositions of Book VII,
the number theoretical part of Theaetetus' proof has no further difficulties.
One can, for instance, reason as follows:
$q$ and $p$ in (1) can of course be taken to be relatively prime.
Then it follows from VII 27 that $p^2$ and $q^2$ (and $p^3$ and $q^3$ as well) are relatively prime,
so that, by VII 21, they are the least of the numbers in that ratio.
But in the left member of (1) occur $n$ and 1 which are also relatively prime and
hence also the least in this ratio.
It follows that $n = p^2$ and $1 = q^2$, so that $n$ is a square.
The case of cubes proceeds analogously.
In my opinion the merit of Theaetetus lies therefore
not in his contribution to the theory of numbers,
but in his study of incommensurable line segments which produce commensurable squares. p.168

**9.3.2.** *Reconstruction by Knorr, 1975*

THEOREM 5. If a square contains an integral number of unit areas, and
if that number is a square integer,
then the side of the square will be commensurable with the unit-line;
but if that number is non-square,
then the side will be incommensurable with the unit. *(Theaetetus,* 148A)
[*Theorem* 5. If a, b are lines and N non-square number such that $a^2=Nb^2$,
then a, b are incommensurable]
…
THEOREM 7. Two lines are commensurable in length if and only if
the squares constructed on them have the ratio which a square integer has to a square integer. *(Elements* X,9)
[*Theorem* 7. If a,b are lines, and A,C are natural numbers, such that AC is not a square number, and $Aa^2=Cb^2$,
then a, b are incommensurable] p.225

The theorems admit of completely straightforward demonstrations,
as soon as the relevant general theorems on relative primes are available,
as we shall show.
Portions of that theory were already implicit in Archytas' proof of the theorem on epimoric ratios, namely, the examination of the techniques of reducing the terms of
ratios.
Hence, by the first or second decade of the fourth century,
the number theory from which the proofs of Theorems 5 and 6 follow as easy deductions
was in preparation. We are thus justified in accepting the most obvious interpretation of
Plato's account: that Theaetetus himself instituted the proofs of the theorems which he is
made to enunciate in the dialogue. p.226

**9.3.3.** *A rendering of the proof of the theorem of Theaetetus on quadratic incommensurabilities suggested by Zeuthen, van der Waerden, and Knorr*

*Proposition 9.3.1.* If $A \cdot a^2 = C \cdot b^2$, and $A \cdot C$ is not a square number, then a, b are incommensurable.

*Proof.* Suppose a, b are commensurable.
Then $a = m \cdot c$, $b = n \cdot c$ for some natural numbers m, n, and line c.
We may assume by the arithmetical Proposition VII.21 that m, n are relatively prime.
Then $A \cdot m^2 = C \cdot n^2$.
By Proposition VII.19, $m^2 / n^2 = C / A$.
By Proposition VII.27, $m^2$, $n^2$ are relatively prime.
By Proposition VII.20, there is a natural number k, such that $C = k \cdot m^2$, $A = k \cdot n^2$.
Then $A \cdot C = k^2 \cdot m^2 \cdot n^2$, a square number, contradiction.

**9.4.** *The strict neutrality thesis of strict Burnyeat on the method of proof employed by Theaetetus*

Although most historians of Greek Mathematics lean towards an arithmetical method of proof for Theaetetus' incommensurabilities, nevertheless Burnyeat, 1978, rejected all reconstructions proposed till then [1978], whether arithmetical or anthyphairetic, on the ground that Plato's *Theaetetus* passage, our only ancient source for Theaetetus' discovery, gives no information on the method followed, either by Theodorus or Theaetetus; in fact, Burnyeat expressed the view that

> Plato has no motive to indicate to the reader whether he has in mind any definite method of proof. Burnyeat, 1978, p.505.

Indeed *Burnyeat, 1978* writes:

> In the face of this impasse the question that needs to be asked is whether Plato has any reason to leave a hint, let alone so indeterminate and ambiguous a hint, as to the mathematical methods used in Theodorus' lesson. As every reader of the dialogue knows, the mathematical scene illustrates a point about definition and examples. When Theaetetus is first asked what knowledge is, he replies by giving examples of knowledge: geometry and the other mathematical sciences he is learning with Theodorus, cobblery and other crafts-each and all of these are knowledge (146cd). Socrates puts him right with an analogy: his answer is like that of someone who, on being asked what clay is, replies, "There is potters' clay, brickmakers' clay, and so on, each and all of which are clay,"

> giving a list of clays instead of the simple, straightforward answer, "It is earth mixed with liquid" (146d-147c). It is at this point that Theaetetus says, "It looks easy now, Socrates, when you put it like that" (147c), and proceeds to tell his story. Theodorus' part in the story does not depend on whether or not he could continue past 17; his role is to provide examples of incommensurability. His case-by-case proof of their incommensurability is mentioned because, if one is not going by a general definition or rule of the kind the boys devised, it is only via construction and proof that examples of incommensurability are forthcoming: construction to obtain a length such as $\sqrt{3}$, which is not marked on any ruler, and proof to show that, divide how you will, you can find no unit to measure without remainder both it and a 1-foot line.
> Beyond that, *Plato has no motive to indicate to the reader whether he has in mind any definite method of proof* or any particular cause for Theodorus to stop at 17. This is not to deny, of course, the legitimacy of speculating about what methods would be available to Theodorus or other fifth-century mathematicians for proving various cases of incommensurability. But this must be an independent inquiry; there is no good reason to expect that the answer is to be squeezed out of one ambiguous sentence in Plato's dialogue. p.505

We will discuss Burnyeat's neutrality thesis in Section 9.6, below.

## 9.5. *Critique of the translations by van der Waerden and Knorr*

These translations are problematical in at least two counts.
The first difficulty is with the almost universal rendering

> since the number of roots appeared to be infinite (9.2.1, van der Waerden)
> since we recognized the powers to be unlimited in number (9.2.2, Knorr)
> since the powers seemed to be unlimited in number (9.2.3, Burnyeat)

of the sentence

> *ἐπειδὴ ἄπειροι τὸ πλῆθος αἱ δυνάμεις ἐφαίνοντο,* 147d7-8

In Negrepontis, 2018 and Negrepontis, Farmaki, Brokou, 2021, arguments have been presented, according to which the only conceivable way to make sense of the sentence 147d7-8 is to read the plural "the powers" in a *distributive* sense:

> since *each power* was seen to be infinite in multitude

It is then strongly suggested that *each power* was seen to be *infinite* in multitude, as a result of Theodorus' lesson, because Theodorus was subjecting each power to anthyphairetic division, and thus proving, by means of Proposition X.2 of the *Elements,* that each power, of those he examined, was incommensurable to the one foot line.

The second difficulty is the almost universal trivialized rendering to try to collect them under one name,
by which we could henceforth call all the roots. (9.2.1, van der Waerden)
we might try to collect them under a single name,
by which we would designate all these powers. (9.2.2, Knorr)
to try to collect the powers under one term
by which we could refer to them all. (9.2.3, Burnyeat)

of the sentence

*πειραθῆναι συλλαβεῖν εἰς ἕν,*
*ὅτῳ πάσας ταύτας προσαγορεύσομεν τὰς δυνάμεις.* 147d8-e1

In the two papers quoted above, the corresponding distributive plural gives the following rendering:

to try to collect [the infinite multitude] into one (συλλαβεῖν εἰς ἕν),
by which (ὅτῳ) we could refer to *each of all these powers* (πάσας ταύτας τὰς δυνάμεις).

Thus Theaetetus' mathematical discovery consists in proving that every power, namely every line a such that $a^2 = Nb^2$, for some non-square number N, is divided into an infinite multitude of parts, namely has infinite anthyphairesis with respect to b, hence is incommensurable to b, but in addition all these infinitely many parts can be "collected into one". Here Plato is employing philosophical language to describe a mathematical theorem, in fact his philosophical method of Name and Logos by which he defined the Sophist, is equivalently described as Division and Collection; as we outlined in Subsection 5.3, Platonic Collection into One refers to the fundamental One-and-Many property of intelligible Beings, whose ultimate cause is definitely anthyphairetic periodicity.

**9.6.** *Critique of the arithmetical non-anthyphairetic reconstructions of Theaetetus' theorem described in the Theaetetus 147d3-e1 by van der Waerden and Knorr, and of the neutrality thesis by Burnyeat*

**9.6.1.** *The arithmetical reconstructions of Theaetetus' theorem suggested by van der Waerden and Knorr are unable to explain the relevance of Theaetetus' theorem to Plato's philosophy*

We find difficult to understand how Zeuthen and van der Waerden, on the one hand have suggested that the method of Theodorus was anthyphairetic, employing Proposition X.2 of the *Elements* (cf. van der Waerden 1954, pp. 144-146), and on the other hand have suggested that the method of Theaetetus was arithmetical, employing Proposition VII.27 of the Elements, even though according to Theaetetus 147d3-e1, young Theaetetus had conceived his proof following the lesson by Theodorus. Perhaps Zeuthen, van der Waerden and other historians were misled by the proof of X.9 presented by Euclid in the *Elements*. In fact, Knorr (1975) writes the following:

The essential content of this passage, as it might appear in the survey literature on
the history of mathematics [Footnote 3] is this: Plato distinguishes two theorems -

> that of Theodorus, to the effect that the square roots of 3, 5, ... , 17 are irrational;
> and that of Theaetetus, that the square root of any non-square integer is irrational.
> [. . . ]
> In fact, this is the sense given to it in a scholium to the ninth proposition of the
> tenth book of the Elements ("two lines are commensurable in length if and only if
> the squares generated on them are in the ratio of square numbers.") The scholiast
> calls this "the discovery of Theaetetus." It would appear that his source for this
> claim is precisely our passage from the Theaetetus. For he cites it and interprets it
> as asserting merely a restricted form of the Euclidean theorem. Now, accepting this
> natural reading of the passage,
> [. . . ] pp. 64-65
>
> [Footnote 3] Accounts of this basic type appear in the histories by W. W. R. Ball
> (1960), J. Tropfke (1921-1924), J. E. Hofmann (1953-57), F. Lasserre (1966), and C. Boyer
> (1968) p. 97

Regardless of the reason which led to van der Waerden's misconception, and despite the fact that he did realize the logical errors and gaps in Euclid's proof, he did not consider that the original method of proof by Theaetetus could have been anything other than arithmetical. Knorr (1975) is at least self-consistent in adopting throughout, for the Pythagoreans, Theodorus and Theaetetus, arithmetical methods. Van der Waerden (1954) writes:

> In the dialogue Plato gives an example of a mathematical investigation made by
> Theaetetus. The example is rather dragged in; it has to serve as an introduction to
> a philosophical discussion, but it does not fit very well. p. 166

It is clear from this statement that van der Waerden does not realize the close connection between Theaetetus' mathematical discovery and Plato's philosophy. The theorem of Theaetetus, reported in the *Theaetetus* 147d-e, does establish the incommensurability of quadratic powers, but it accomplishes not only this but something much more, it provides complete knowledge of the power; by proving anthyphairetic periodicity, by Proposition X.2, we immediately have incommensurability, but by obtaining at a finite stage the period of the anthyphairesis we also have complete knowledge of an infinite entity at a finite stage. Thus the mathematical example is not "rather dragged", and it is not the case that "it does not fit very well", but on the contrary it reveals the anthyphairetic basis of Plato's philosophy. The same remarks apply for Knorr's somewhat different, but also arithmetical and non-anthyphairetic reconstruction.

**9.6.2.** *The neutrality thesis of Burnyeat fails to take into account (a) the passage* Theaetetus *148c-d on the urge to imitate Theaetetus' method, indicating great interest by Plato of Theaetetus' method, and (b) the anthyphairetic nature of the actual imitation, as outlined in Sections 5-7*

Burnyeat's translation of the *Theaetetus*' 147d-e (given in 9.2.3) is almost identical to the translation by van der Waerden, Knorr, and is thus subject to the same critique that theirs is. Burnyeat's thesis (presented in 9.4) that Plato has no motive to provide any information of Theaetetus' method is most certainly wrong. For how would we reconcile Plato's lack of motive or interest to provide some hint of the method of proof by Theaetetus with Plato's avowed urge and intention to imitate this method of Theaetetus on the subject of the Knowledge of the intelligible Beings (*Theaetetus* 148d4-5), a subject of supreme importance for Plato's philosophy? Moreover, Burnyeat's neutrality thesis (9.4) fails to take into account the anthyphairetic nature of the imitation, realized in the *Sophist* and the *Statesman*.

## 9.7. *The problem with the philosophical interpretations of Plato's* Logos

The most fundamental concept in Plato's philosophy is that of "Logos". As Plato writes in the *Sophist* 260a-b

> if we were deprived of this [Logos], we would be deprived of philosophy, which would be the greatest calamity

As we have indicated in Section 5, employing passages from the *Sophist* and the *Meno*, the standard interpretations of Logos, as "account", "discourse", "rational definition", "rational explanation" and "reason" are not the correct ones. Plato's Logos is an imitation of Theaetetus' Logos Criterion for the periodicity of a geometric anthyphairesis, providing not only incommensurability, but complete Knowledge of its infinite structure and Equalization of all its parts, collecting all in a One, and thus making it an ideal candidate for an intelligible Being.

*Theaetetus* 202c2-5

> for he who cannot give and receive Logos (logon) of a thing is without knowledge (anepistemona) of it; but when he has acquired Logos (logon) he may possibly have become all that I have said and may now be perfect in knowledge (epistemen).

*Note*. Logos is rendered in the classical translation by Fowler (1921) as rational explanation.

*Sophist* 260a5-b2, *Logos* is necessary for the Knowledge of intelligible Beings.

> Stranger Our object was to establish that the Logos (logon) is one of the kinds of the Beings (ton onton). For if we were deprived of this [Logos], we would be deprived of philosophy, which would be the greatest calamity; moreover, we must at the present come to an agreement (diomologesasthai) about the nature (ti pot' estin) of Logos (λόγον), and if we were robbed of it [Logos], there would be no true Being (med' einai) and we could no longer be able to provide Logos (legein), and we should be robbed of it [Logos], if we agreed that there is no mixture (meixin) of anything with anything.

*Note*. The classical translation by Fowler (1921) renders all six references to "Logos" as "discourse".

*Statesman* 286a4-7, Logos is necessary for the Knowledge of intelligible Beings.

> We must therefore study so as to have the power (dunaton einai) to give and receive (dounai kai dexasthai) the Logos (logon) of each of them; for the immaterial [intelligible] Beings, which are the noblest and greatest, can be made known (saphos deiknutai) by Logos (logoi) only, and it is for their [the intelligibles'] sake (heneka) that all Logoi are provided (legomena).

*Note*. The classical translation by Fowler and Lamb (1925) renders "logon" as "rational definition", "logoi" as "reason", and "legomena" as "all we are saying".

*Republic* 534b3-6

"And do you not also give the name dialectician (dialectikon) to the man who is able to exact an account (logon) of the true Being of each thing? And will you not say that the one who is unable to do this, in so far as he is incapable of rendering Logos (logon) to himself and others, does not possess intelligence about the matter?"

*Note*. Logos is rendered in the classical translation by Shorey (1935) as "account".

The renderings of the crucial term "logos" in the *Meno* 86b7 and 86c2 by Platonic scholars as "argument" and "word", respectively.

> There are other things about the argument that I would not confidently affirm, but that we shall be better men, more courageous and less idle, if we think we ought to inquire into what we do not know, instead of thinking that because we cannot find what we do not know we ought not to seek it – that I would do battle for, so far as possible, in word and deed. [86b6-c2, in Halper (2007), p.239]

> For he goes on asserting, forcefully and without any ambiguity, that in all that was said there was only one point (81 d 5 - e 1; 86 b 2 - 4) for which he would fight strenuously, in word and deed (kai logōi kai ergoi), whenever he was able to do so. [Klein (1965), p. 182]

> As for the other points, I wouldn't absolutely insist on the argument. But I would fight, both in word and deed, for the following point: that we would be better, more manly, and less lazy if we believed that we ought to inquire into what we do not know, than if we believed that we cannot discover what we do not know and so have no duty to inquire. [86b6–c2, in Scott (2005), p. 121]

Logos is rendered as argument and word.

> socrates: Yes, I think so too, Meno. I wouldn't support every aspect of the argument with particular vigour, but there's one proposition that I'd defend to the death, if I could, by argument and by action: that as long as we think we should search for what we don't know we'll be better people– less fainthearted and less lazy– than if we were to think that we had no chance of discovering what we don't know and that there's no point in even searching for it. [Waterfield (2005), p. 124]

Logos is rendered as argument.

> So I seem to myself, Meno. All the other points I have made in support of the argument are not such as I can confidently assert; but that, if we are convinced we should inquire after what we do not know, we should get better and braver and less lazy than if we believe that we are neither able nor obliged to inquire after things we do not know – this is something for which I am determined to fight, so far as I am able, both in word and deed. [86b6-c2, in Ebert (2007), p. 197]

> socrates: I think so too, Meno. In defense of the argument I would not affirm the other points very strongly, but that we would be better, more manly, and less lazy by believing that one should search for what one doesn't know than if we believed that we cannot discover what we do not know and should not even search for it – that is something over which I would fiercely contend, if I were able, in both word and deed. [Sedley and Long, (2011), p. 23]

> As far as the other points (ta alla) are concerned, I wouldn't altogether take a stand on the argument. But that we will be better and more manly and less idle if we think one should inquire into that which one doesn't know than if we think it isn't possible to discover what we don't know and that we don't need to inquire into it: that is something I would certainly fight for to the end, if I could, in both word and deed. [85d12–86c2, in Fine (2014), p. 150]

> I would not confidently assert the other things said in defense of this account, but that we would be better and braver and less idle if we believe that one ought to inquire concerning those things he fails to know than if one believes it is not possible to discover nor necessary to inquire concerning those things one fails to know (εἰ οἰοίμεθα ἃ μὴ ἐπιστάμεθα μηδὲ δυνατὸν εἶναι εὑρεῖν μηδὲ δεῖν ζητεῖν), I would fight for in both word and deed as far as I am able. (Meno 86B6–C2; emphasis added) [Benson (2015), p. 62]

Nehamas (1984, p.36) expresses the opinion that the meaning that Plato has in mind with Logos must come from a study of the Divisions in the *Sophist* and the *Statesman*.

> The answer [to the meaning of Logos in Plato] can only be given using the long and complicated divisions through which the Sophist and the Politicus attempt to define their subject matter. Both dialogues try to define a species of art (techne) or science (episteme), that of the sophist and that of the statesman (cf. Pol, 258b2-7). In both dialogues, the art in question is defined by being located within a determinate network of other arts and sciences. And the process by means of which the definition is reached provides, at the same time, both an explanation of the aporetic ending of the Theaetetus and an Illustration of the obscure opening of the Philebus. Plato, therefore, does try to answer the question of Theaetetus. But his answer is not cryptically contained within the dialogue itself, either negatively (as Cornford argued) or positively (as Fine suggested). His answer is given in the two dialogues that follow the Theaetetus.

**Section 10.** *The reconstruction of Theaetetus' theory of ratios of magnitudes only for ratios with no recourse to Eudoxus' condition (Definition 4 in Book V of the* Elements*)*

The detailed reconstruction of the proof of Theaetetus' theorem reveals that Theaetetus' theory of *magnitudes* is essentially a theory of ratios of *lines*, designed principally for the needs of periodic anthyphairesis. Thus the Logos Criterion for anthyphairetic periodicity is essentially just a restatement in more suggestive and convenient language of the more heavy and burdensome Pythagorean preservation of Application of Areas/Gnomons. Thus, Theaetetus' theory is for the limited class of those ratios of lines with *finite or eventually periodic anthyphairesis*. Exactly because of the limited nature of the theory, it is possible to prove, with the help of the crucial Proposition 8.1.4, all the basic properties of this theory without recourse to Eudoxus' condition. The extension of the theory to areas and surfaces is achieved with the help of the analogue of the key Proposition VI.1 of the *Elements*. The fact that this theory is really a theory not for general magnitudes, but principally for lines, makes necessary the many special cases in the proof of its Propositions. This feature is pointed out by Aristotle in *Analytics Posterior* 74a, and the vehicle of passing from lines to areas and surfaces is exactly the proposition (analogue to Proposition VI.1 of the *Elements*) singled out by Aristotle in the *Topics* 158b to inform us about this theory.

## 10.1. *The introduction of the generalized side and diameter numbers*

**10.1.1.** *Definition. The generalized side-diameter numbers $p_n$, $q_n$*
Let there be an infinite anthyphairesis of a to b:
$a = k_0 \cdot b + e_1, \quad e_1 < b$
…
$e_{m-2} = k_{m-1} \cdot e_{m-1} + e_m, \quad e_m < e_{m-1}$
$e_{m-1} = k_m \cdot e_m + e_{m+1}, \quad e_m < e_{m+1}$
$e_m = k_{m+1} \cdot e_{m+1} + e_{m+2}, \quad e_{m+1} < e_{m+2}$
…
$e_{n-2} = k_{n-1} \cdot e_{n-1} + e_n, \quad e_n < e_{n-1}$
$e_{n-1} = k_n \cdot e_n + e_{n+1}, \quad e_{n+1} < e_n$
… .
We set
$p_1 = 1$, $q_1 = k_0$, and
$p_n = k_{n-1} \cdot p_{n-1} + p_{n-2}$, $q_n = k_{n-1} \cdot q_{n-1} + q_{n-2}$ for every n.

*Note.* The (classical) side and diameter numbers are identified as the generalized side and diameter numbers of the anthyphairesis of the diagonal (diameter) to side of a square.
[Indeed the anthyphairesis of diameter to side of a square is
$k_0 = 1$, $k_n = 2$ for all n = 1, 2, … .
Then the generalized side-diameter numbers are recursively defined by
$p_1 = 1$, $q_1 = k_0 = 1$, and,
$p_{n+2} = 2 \cdot p_{n+1} + p_n$, $q_{n+2} = 2 \cdot q_{n+1} + q_n$.
Then
$p_{n+2} = 2 \cdot p_{n+1} + p_n = p_{n+1} + p_{n+1} + p_n = p_{n+1} + (p_n + q_n) + p_n = p_{n+1} + q_{n+1}$, and
$q_{n+2} = 2 \cdot q_{n+1} + q_n = q_{n+1} + q_{n+1} + q_n = q_{n+1} + 2 \cdot p_n + q_n + q_n = q_{n+1} + 2 \cdot (p_n + q_n) = 2 \cdot p_{n+1} + q_{n+1}$,
and the recursive definition of the side and diameter numbers is derived
$p_{n+2} = p_{n+1} + q_{n+1}$,
$q_{n+2} = 2 \cdot p_{n+1} + q_{n+1}$].

A major difference between the classical side and diameter numbers and the generalized ones is that the recursive rule defining the generalized side and diameter numbers is *far inferior* to the classical rule, in that at the stage n they are defined in terms of the *previous two stages*, n-1 and n-2, and not just in terms of the immediately previous one n-1.

**10.1.2.** *The generalized side and diameter numbers are useful in expressing the anthyphairetic remainders in terms of the initial magnitudes a and b.*

*Proposition (the anthyphairetic remainders in terms of the generalized side and diameter numbers).* With the notation of Definition 10.1.1, we have
$e_n = (-1)^n (q_n \cdot b - p_n \cdot a)$ for all n = 1, 2, … .

*Proof.* $e_1 = a - k_0 \cdot b = (-1)^1 (q_1 \cdot b - p_1 \cdot a)$.
Assume inductively that:
$e_{n-2} = (-1)^{n-2} (q_{n-2} \cdot b - p_{n-2} \cdot a)$, and $e_{n-1} = (-1)^{n-1} (q_{n-1} \cdot b - p_{n-1} \cdot a)$.
Then:
$e_n =$
$e_{n-2} - k_{n-1} \cdot e_{n-1} =$
$(-1)^{n-2} (q_{n-2} \cdot b - p_{n-2} \cdot a) - (-1)^{n-1} \cdot k_{n-1} \cdot (q_{n-1} \cdot b - p_{n-1} \cdot a) =$
$(-1)^n [q_{n-2} \cdot b - p_{n-2} \cdot a + k_{n-1} \cdot (q_{n-1} \cdot b - p_{n-1} \cdot a)] =$
$(-1)^n [(q_{n-2} + k_{n-1} \cdot q_{n-1}) \cdot b - (p_{n-2} + k_{n-1} \cdot p_{n-1}) \cdot a] =$
$(-1)^n (q_n \cdot b - p_n \cdot a)$.

**10.1.3.** *Proposition.* If $A \cdot a^2 = B \cdot a \cdot b + C \cdot b^2$, and $a \cdot d = b \cdot c$, then $A \cdot c^2 = B \cdot c \cdot d + C \cdot d^2$.

*Proof.* From the assumption we have $A \cdot a^2 = B \cdot a \cdot b + C \cdot b^2$,
hence $A \cdot a^2 \cdot c = B \cdot a \cdot b \cdot c + C \cdot b^2 \cdot c$
hence $A \cdot a^2 \cdot c = B \cdot a \cdot b \cdot c + C \cdot b \cdot a \cdot d$,
hence $A \cdot a \cdot c = B \cdot b \cdot c + C \cdot b \cdot d$,
hence $A \cdot a \cdot c^2 = B \cdot b \cdot c^2 + C \cdot b \cdot d \cdot c$,
hence $A \cdot a \cdot c^2 = B \cdot c \cdot a \cdot d + C \cdot a \cdot d^2$,
hence $A \cdot c^2 = B \cdot c \cdot d + C \cdot d^2$.

**10.1.4.** *Proposition.*
[1] For every finite ordered sequence of natural numbers $k_0, k_1, …, k_n$,
there are lines a and b, a > b, such that
    Anth(a, b) = [period ($k_0, k_1, …, k_n$)], and
    $p_{n+1} \cdot a^2 = (q_{n+1} - p_n) \cdot a \cdot b + q_n \cdot b^2$.
[2] If lines c and d satisfy Anth(c, d) = [period ($k_0, k_1, …, k_n$)], then
    $p_{n+1} \cdot c^2 = (q_{n+1} - p_n) \cdot c \cdot d + q_n \cdot d^2$.

*Proof.* [1] The first n + 1 steps of the anthyphairesis of a to b are:
    a  = $k_0 \cdot b + e_1$,     $e_1 < b$
    b  = $k_1 \cdot e_1 + e_2$,     $e_2 < e_1$
    $e_1$ = $k_2 \cdot e_2 + e_3$,     $e_3 < e_2$
    …

$$e_{n-2} = k_{n-1} \cdot e_{n-1} + e_n, \quad e_n < e_{n-1}$$
$$e_{n-1} = k_n \cdot e_n + e_{n+1}, \quad e_{n+1} < e_n.$$

We also impose the cross product condition
$$a \cdot e_{n+1} = b \cdot e_n.$$
The cross product condition implies the periodicity of the anthyphairesis; indeed, by Proposition 7.1.4. $Anth(a, b) = Anth(e_n, e_{n+1})$, and thus by induction:

$$Anth(a, b) = [k_0, k_1, \ldots, k_n, Anth(e_n, e_{n+1})]$$
$$= [k_0, k_1, \ldots, k_n, Anth(a, b)]$$
$$= [k_0, k_1, \ldots, k_n, k_0, k_1, \ldots, k_n, Anth(e_n, e_{n+1})]$$
$$= \ldots$$
$$= [\text{period }(k_0, k_1, \ldots, k_n)].$$

The cross product relation further insures that a, b satisfy an arithmetical quadratic equality in excess; indeed, by Proposition 10.1.2, $e_n = (-1)^n (q_n \cdot b - p_n \cdot a)$, and $e_{n+1} = (-1)^{n+1}(q_{n+1} \cdot b - p_{n+1} \cdot a)$, Hence, $a \cdot e_{n+1} = (-1)^{n+1} \cdot a \cdot (q_{n+1} \cdot b - p_{n+1} \cdot a)$, and $b \cdot e_n = (-1)^n \cdot b \cdot (q_n \cdot b - p_n \cdot a) = (-1)^{n+1} b(p_n \cdot a - q_n \cdot b)$. Thus, the cross product condition takes the form: $a \cdot (q_{n+1} \cdot b - p_{n+1} \cdot a) = b \cdot (p_n \cdot a - q_n \cdot b)$, hence
$$p_{n+1} \cdot a^2 = (q_{n+1} - p_n) \cdot a \cdot b + q_n \cdot b^2.$$

[2] Since $Anth(a, b) = Anth(c, d)$, we have that the anthyphairesis of c to d is purely periodic. Since the generalized side and diameter numbers depend solely on the sequence of quotients $k_0, k_1, \ldots, k_n$ it follows that the generalized side $r_i$ and diameter $s_i$ numbers of c to d coincide with the generalized side $p_i$ and diameter $q_i$ numbers of a to b, respectively, namely $r_i = p_i$, $s_i = q_i$ for i = 1, 2, ..., n, n+1;

hence, denoting by $f_i$ the $i^{th}$ anthyphairetic remainder of c to d, it follows that: $c \cdot f_{n+1} = d \cdot f_n$, hence, by Proposition 10.1.3, c, d satisfy the same arithmetical quadratic equation
$p_{n+1} \cdot c^2 = (q_{n+1} - p_n) \cdot c \cdot d + q_n \cdot d^2$.

From this result by completing the anthyphairesis of both a to b and c to d upwards we obtain, without difficulty, the general case. The details are left to the reader.

## *10.2. Theaetetus' theory of ratios of magnitudes*

Proposition 10.1.4 makes possible the proof of the fundamental Proposition 10.2.4: in Theaetetus' theory, for lines a, b, c, d, the condition of proportion a/b = c/d is equivalent to the equality of the cross products $a \cdot d = b \cdot c$, with no recourse to Eudoxus' condition. In turn, following essentially the old idea of Becker, 1933, we can prove not only the *Alternando* property for lines, but also the crucial analogue of Proposition V.8 of the *Elements*.

**10.2.1.** *Definition of Theatetean ratio.*
A line a possesses Theaetetean ratio with respect to line b if the anthyphairesis of a to b is finite or eventually periodic.

**10.2.2.** *Definition of proportion.*
If lines a, b and c, d have finite or eventually periodic anthyphairesis, then the ratios a/b and c/d are *analogous (in proportion)* a/b = c/d in the Theaetetean sense if $Anth(a, b) = Anth(c, d)$.

**10.2.3.** *Proposition* (analogue of the Transitive property, proposition V.11 for lines, no analogue Proposition for numbers).
Let a, b and c, d and e, f be pairs of lines possessing Theaetetean ratios, namely with finite or eventually periodic anthyphairesis. Then if a/b = c/d, and c/d = e/f, then a/b = e/f;

*Proof.* By definition 10.2.1, either Anth(a, b) = [$k_0, k_1, \ldots, k_n$, period ($k_{n+1}, k_{n+2}, \ldots, k_m$)] if it is eventually periodic, or Anth(a, b) = [$k_0, k_1, \ldots, k_n$] if it is finite.
Also, by definition 10.2.2, a/b = c/d means that Anth(a, b) = Anth(c, d).
And, by definition 10.2.2 again, c/d = e/f means that Anth(c, d) = Anth(e, f).
Thus Anth(a, b) = Anth(e, f) and this, by the same definition, means that a/b = e/f.

**10.2.4.** *Fundamental Proposition* (analogue of Proposition VI.16 for lines, Proposition VII.19 for numbers).
Let a, b and c, d be pairs of lines possessing Theaetetean ratios, namely with finite or eventually periodic anthyphairesis.
Then a/b = c/d (namely Anth(a, b) = Anth(c, d)) if and only if a·d = b·c.

*Proof.* Let a, b and c, d be two pairs of lines with finite or eventually periodic anthyphairesis.
If a·d = b·c, then by Proposition 5.2.1, Anth(a, b) = Anth(c, d), namely a/b = c/d.
If a/b = c/d, then by Proposition 10.1.4, a, b and c, d satisfy the same arithmetical quadratic equation, say $A \cdot a^2 = B \cdot a \cdot b + C \cdot b^2$ and $A \cdot c^2 = B \cdot c \cdot d + C \cdot d^2$. By Proposition 7.1.2, a·d = b·c.

**10.2.5.** *Proposition* (analogue of Proposition V.9 for lines, no analogue Proposition for numbers is explicitlystated in Book VII but it is implicitly used).
If lines a, b and a, c have finite or eventually periodic anthyphairesis, and a/b = a/c, then b = c.

*Proof.* By definition a/b = a/c means that Anth(a, b) = Anth(a, c).
By Proposition 10.2.4, a·c = a·b, hence b = c.

*Note*. For the restrictive class of lines a, b and c, d with *finite or eventually periodic anthyphairesis*, the proof of the Fundamental Propositions 10.2.4 and 10.2.5 are proved *without* any use of Eudoxus' condition Definition V.4 of the *Elements,* thus avoiding the anachronism that plagues the earlier reconstructions. This by itself is a strong argument in favor of our reconstruction.

Regarding the next proposition we will prove, the *Alternando*, van der Waerden (1975, p. 177) writes:

> But one property, viz. the interchange of the means, causes difficulty. And now it is curious that, according to Aristotle, it was indeed exactly the proof of this proposition, which at first led to difficulties. "Formerly", says Aristotle in Anal. Post. 15, "this proposition was proved separately for numbers, for line segments, for solids and for periods of time. But after the introduction of the general concept which includes numbers as well as lines, solids and periods of time" (viz. the concept of magnitude), "the proposition could be proved in general".

The proof for numbers can be found in Book VII (VII 13). What was the old proof for line segments?

O. Becker has advanced an ingenious hypothesis for this. From the proportionality
(1) a:b = c:d,
one deduces first the equality of the areas
(2) ad=bc,
then interchanges b and c, and finally returns to the proportionality
(3) a:c = b:d.

Van der Waerden refers to the following remarkable Aristotle's *Analytics Posterior* 74a17-25 passage:

> καὶ τὸ ἀνάλογον ὅτι καὶ ἐναλλάξ,
> ᾗ ἀριθμοὶ καὶ ᾗ γραμμαὶ καὶ ᾗ στερεὰ καὶ ᾗ χρόνοι,
> ὥσπερ ἐδείκνυτό ποτε χωρίς,
> ἐνδεχόμενόν γε κατὰ πάντων μιᾷ ἀποδείξει δειχθῆναι·
> ἀλλὰ διὰ τὸ μὴ εἶναι ὠνομασμένον τι ταῦτα πάντα ἓν,
> ἀριθμοί μήκη χρόνοι στερεά, καὶ εἴδει διαφέρειν ἀλλήλων, χωρὶς ἐλαμβάνετο.
> νῦν δὲ καθόλου δείκνυται·
> οὐ γὰρ ᾗ γραμμαὶ ἢ ᾗ ἀριθμοὶ ὑπῆρχεν, ἀλλ' ᾗ τοδί,
> ὃ καθόλου ὑποτίθενται ὑπάρχειν.

> Again, the law that proportionals alternate might be supposed to apply
> to numbers qua numbers, and similarly to lines, solids and periods of time;
> as indeed it used to be demonstrated of these subjects separately.
> It could, of course, have been proved of them all by a single demonstration,
> but since there was no single term to denote the common quality of
> numbers, lengths, time and solids,
> and they differ in species from one another, they were treated separately;
> but now the law is proved universally;
> for the property did not belong to them
> qua lines or qua numbers,
> but qua possessing this special quality which they are assumed to possess universally.
> [translation by Tredennick ()]

We will follow Becker's idea for the proof of proposition

**10.2.6.** *Proposition (Alternando,* analogue of Proposition V.16 for lines, Proposition VII.13 for numbers).
If a/b = c/d and Anth(a, c), Anth(b, d) are finite or eventually periodic,
then a/c = b/d.

*Proof.* By Proposition 10.2.4, a·d = b·c. If Anth(a, c) is finite or eventually periodic, then by Proposition 10.2.4, a/c = b/d.

**10.2.7.** *Proposition (Ex Equali,* analogue of Proposition V.22 for lines, Proposition VII.14 for numbers).
If a, b, c, d, e, f are lines and a/b = d/e and b/c = e/f and Anth(a, c), Anth(d, f) are finite or eventually periodic, then a/c = d/f.

*Proof.* By the Fundamental Proposition (10.2.4), we have:
a·e = b·d, and
b·f = c·e.
Then a·b·f = a·c·e = b·c·d,
hence a·f = c·d.
By the Fundamental Proposition (10.2.4), a/c = d/f.

**10.2.8.** *Proposition* (perturbed proportion, analogue of Proposition V.23 for lines).
If a, b, c, d, e, f are lines and a/b = e/f and b/c = d/e, and Anth(a, c), Anth(d, f) are finite or eventually periodic, then a/c = d/f.

*Proof.* By the Fundamental Proposition (10.2.4), we have
a·f = b·e and b·e = c·d.
By transitivity (10.2.3), a·f = c·d.
By the Fundamental Proposition (10.2.4), we have
a/c = d/f.

**10.2.9.** *Proposition* (part [1] is the analogue of Proposition V.12 for lines and Proposition VII.12 for numbers while part [2] is the analogue of Proposition V.19 for lines and Proposition VII.11 for numbers).
If a, b, c, d are lines and Anth(a, b) = Anth(c, d) is finite or eventually periodic anthyphairesis, then:
[1] Anth(a + c, b + d) is finite or eventually periodic, and (a + c)/(b + d) = a/b, and
[2] if a > c and b > d, then Anth(a – c, b – d) is finite or eventually periodic, and (a – c)/(b – d) = a/b.

*Proof.*
[1] By proposition 10.2.4, and since Anth(a, b) = Anth(c, d), we have: a·d = b·c.
Then: a·d + a·b = b·c + a·b,
a·(b + d) = b·(a + c).
By proposition 5.2.1, we have that Anth(a, b) = Anth(a + c, b + d) and then, by definition 10.2.2 (a + c)/(b + d) = a/b.

[2] By proposition 10.2.4, and since Anth(a, b) = Anth(c, d), we have: a·d = b·c.
Then: a·b + a·d = a·b + b·c,
a·b – b·c = a·b – a·d,
b·(a – c) = a·(b – d).
By proposition 5.2.1, we have that Anth(a, b) = Anth(a – c, b – d) and then, by definition 10.2.2 (a – c)/(b – d) = a/b.

**10.2.10.** *Proposition* (part [1] is the analogue of Proposition V.18 for lines and part [2] is the analogue of Proposition V.17 for lines, but there is no analogue for numbers).
If a, b, c, d are lines and Anth(a, b) = Anth(c, d) is finite or eventually periodic anthyphairesis, then:
[1] Anth(a + b, b) = Anth(c + d, d) is finite or eventually periodic, and (a + b)/b = (c + d)/d, and
[2] if a – b > b and c – d > d, then Anth(a – b, b) = Anth(c – d, d) is finite or eventually periodic, and (a – b)/b = (c – d)/d.

*Proof.*
[1] Consider the case where Anth(a, b) = Anth(c, d) and both are finite.
Let their anthyphairesis be denoted by $[k_0, k_1, \ldots, k_n]$.
Then, it is easy to demonstrate, that the anthyphairesis of a + b to b and of c + d to d, differs from the anthyphairesis of a to b and c to d, respectively, only in the first quotient.
Essentially, instead of $k_0$ the first quotient now becomes $k_0 + 1$.
Then Anth(a + b, b) = Anth(c + d, d) = $[k_0 + 1, k_1, \ldots, k_n]$ and, by definition 10.2.2,
(a + b)/b = (c + d)/d
Similarly, we can demonstrate that the proposition is true when the anthyphairesis of a to b is periodic or eventually periodic.

[2] Consider the case where Anth(a, b) = Anth(c, d) and both are finite.
Let their anthyphairesis be denoted by $[k_0, k_1, \ldots, k_n]$.
Then, it is easy to demonstrate, that the anthyphairesis of a – b to b and of c – d to d, differs from the anthyphairesis of a to b and c to d, respectively, only in the first quotient.
Essentially, instead of $k_0$ the first quotient now becomes $k_0 – 1$.
Then Anth(a + b, b) = Anth(c + d, d) = $[k_0 – 1, k_1, \ldots, k_n]$ and, by definition 10.2.2,
(a – b)/b = (c – d)/d
Similarly, we can demonstrate that the proposition is true when the anthyphairesis of a to b is periodic or eventually periodic.

**10.2.11.** *Proposition* (*Topics* Proposition, analogue of VI.1 for lines, Propositions VII.17, 18 for numbers).
If a, b, c are lines, with Anth(a, b) eventually periodic, then a·c/b·c = a/b.

*Proof.* As already noted, we have that Anth(a, b) = Anth(a·c, b·c), hence a·c/b·c = a/b.

**10.2.12.** *Proposition* (analogue of Proposition V.9 for rectilinear areas).
If A, B, C are rectilinear areas, with A/C = B/C, and Anth(A, C) is finite or eventually periodic, then B = C.

*Proof.* By the definition of A/B = A/C, we have Anth(A, B) = Anth(A, C).
Let r be an assumed line.
By Proposition I.45 of the *Elements*, there are lines a, b, c, such that A = a·r, B = b·r, C = c·r.
By Proposition 10.2.11, A/B = a·r/b·r = a/b, A/C = a·r/c·r = a/c.
Hence, by transitivity (10.2.3), a/b = a/c
By the fundamental Proposition 10.2.4, a·c = a·b, hence b = c.
Hence B = b·r = c·r = C.

**10.2.13.** *Proposition* (*Alternando,* analogue of Proposition V.16 for rectilinear areas).
If A, B, C, D are rectilinear areas, with A/B = C/D, and Anth(A, C), Anth(B, D) finite or eventually periodic, then A/C = B/D.

*Proof.* By the definition of A/B = C/D, we have Anth(A,B) = Anth(C, D).
Let r be an assumed line.
By Proposition I.45 of the *Elements*, there are lines a, b, c, d, such that
A = a·r, B = b·r, C = c·r, D = d·r
By Proposition 10.2.11, A/B = a·r/b·r = a/b, C/D = c·r/d·r = c/d.
Hence, by transitivity (10.2.3), a/b = c/d.
By the fundamental Proposition 10.2.4, a·d = b·c.
By Proposition 10.2.6 (*Alternando* for lines) a/c = b/d
By Proposition 10.2.11 and transitivity, A/C = B/D.

**10.2.14.** *Proposition* (Analogue to V.22, part [1] is the *Ex Equali* for rectilinear areas and part [2] is the *Ex Equali* for a mixed proportion).
[1] If A, B, C, D, E, F are rectilinear areas, with A/B = D/E, and B/C = E/F, and Anth(A, C), Anth(D, F) finite or eventually periodic, then A/C = D/F.
[2] If A, B, C are rectilinear areas, d, e, f are lines, with A/B = d/e, B/C = e/f, and Anth(A, C), Anth(d, f) are finite or eventually periodic, then A/C = d/f.

**10.2.15.** *Proposition* (Analogue to V.23, part [1] is the perturbed proportion for rectilinear areas and part [2] is the perturbed proportion for a mixed proportion).
[1] If A, B, C, D, E, F are rectilinear areas, with A/B = E/F, and B/C = D/E, and Anth(A, C), Anth(D, F) finite or eventually periodic, then A/C = D/F.
[2] If A, B, C are rectilinear areas, d, e, f are lines, with A/B = e/f, B/C = d/e, and Anth(A, C), Anth(d, f) are finite or eventually periodic, then A/C = d/f.

*Proofs* of Propositions 10.2.14, 15 are analogous to the proof of Proposition 10.2.12 and are left to the reader.

*Note*. The fragmentation of the proof of a property in separate cases for lines, areas, volumes, as described by Aristotle in the *Analytics Posterior* passage quoted above, is confirmed, not only for the *Alternando* property, but also for the other properties of Theaetetus' theory of proportion.